\documentclass{amsart}
\input BoxedEPS																	% the commands needed
\SetepsfEPSFSpecial				  								%this puts the pictures in for UNIX
\HideDisplacementBoxes										% this hides unsightly boxes in document	

\newtheorem{theorem}{Theorem}[section]
\newtheorem{lemma}[theorem]{Lemma}
\newtheorem{corollary}[theorem]{Corollary}
\newtheorem{conjecture}[theorem]{Conjecture}

\theoremstyle{definition}

\theoremstyle{remark}
\newtheorem{remark}[theorem]{Remark}

\numberwithin{equation}{section}

\newcommand{\HHH}{{\mathbb H}^3}    % hyperbolic 3-space
\newcommand{\HH}{{\mathbb H}^2}     % hyperbolic 2-space
\newcommand{\Zz}{{\mathbb Z}}       % integers
\newcommand{\PSLC}{{\rm PSL_2}({\mathbb C})}  % isometry group PSL_2(C)
\newcommand{\NOE}{N^\circ_\epsilon} % N - {rank one and two cusps}

\begin{document}

\title{Injectivity Radius Bounds in Hyperbolic Convex Cores I}

\author{Carol E.\ Fan}
\address{Department of Mathematics, Oklahoma State University, 
Stillwater, Oklahoma  74078-1058}
\curraddr{Department of Mathematics, Loyola Marymount University, Los Angeles,
CA  90045}
\email{cfan@lmumail.lmu.edu}
\thanks{I would like to thank Richard Canary for advising me during this
project.  Furthermore, Francis Bonahon, Timothy Comar, Mingqing Ouyang,
Peter Scott, and Edward Taylor were very helpful with comments and
discussions.  The results in this paper represent a portion of my Ph.D.\
thesis completed at The University of Michigan.}

%    General info

\date{November 4, 1998}

%\dedicatory{This paper is dedicated to our authors.}


\begin{abstract}
A version of a conjecture of McMullen is as follows:  Given a
hyperbolizable 3-manifold M with incompressible boundary, there exists a
uniform constant K such that if N is a hyperbolic 3-manifold homeomorphic to the
interior of M, then the injectivity radius based at points in the convex core
of N is bounded above by K.  This conjecture suggests that convex cores are
uniformly congested.  In previous
work, the author has proven the conjecture for $I$-bundles over a closed
surface, taking into account the possibility of cusps. In this paper, we
establish the conjecture in the case that M is a book of I-bundles or an
acylindrical, hyperbolizable 3-manifold.  In particular, we show that if M is a
book of I-bundles, then the bound on injectivity radius depends on the number of
generators in the fundamental group of M.    
\end{abstract}

\maketitle

\section{Introduction}

In this paper, we investigate the geometry of convex cores of
hyperbolic 3-manifolds which are homeomorphic to the interior of a book of
$I$-bundles or an acylindrical, hyperbolizable 3-manifold.  Specifically, we
show that if
$M$ is a book of 
$I$-bundles or an  acylindrical, hyperbolizable 3-manifold, then there exists
a uniform upper bound on injectivity radius for points in the convex core
of any hyperbolic 3-manifold homeomorphic to the interior of
$M$.

The main result relies on a theorem of Kerckhoff-Thurston
\cite{kerck/thur} which established the existence of an upper bound on injectivity radius for points
in the convex core of hyperbolic 3-manifolds without cusps where the
manifolds are homotopy equivalent to a closed surface.  The main theorem also
makes use of an extension of their theorem \cite{fan1} which includes the
possibility of cusps.  Our main result is:

\begin{theorem}
Let $M$ be a book of $I$-bundles or an acylindrical, hyperbolizable
3-manifold.  Then there exists a constant $K$ such that if $N$ is a
hyperbolic 3-manifold homeomorphic to the interior of $M$ and $x \in C(N)$,
then $inj_N(x) \leq K$.
\end{theorem} 

The main theorem is also related to a conjecture of McMullen:

\begin{conjecture}
\label{macconj}  (McMullen, \cite{bielefeld}) Let $N$ be a
hyperbolic 3-manifold homotopy equivalent to a compact 3-manifold $M$.  Then
$C(N)$ does not contain an embedded ball of radius $L$, where $L$ depends on
the number of generators of $\pi_1(N)$.
\end{conjecture}

For books of $I$-bundles, we have result which is slightly stronger than
what is predicted by the conjecture, because the injectivity radius measures
the radius of balls embedded in $N$, not in the convex core of $N$:

\begin{corollary} Let $N$ be a hyperbolic 3-manifold homotopy equivalent to a
book of $I$-bundles.  Then there exists a constant $L$ such that for
$x \in C(N)$, $inj_N(x) \leq L$, where $L$ depends on the
number of generators of $\pi_1(N)$.
\end{corollary}
 
The upper bound on injectivity radius, combined with a result of
McMullen's \cite{mcmullen}, can be used to show that the limit set varies
continuously over the space of hyperbolic 3-manifolds homeomorphic to the
interior of a fixed hyperbolizable 3-manifold $M$ under the geometric
topology.  

\begin{corollary}
Let $M$ be a book of I-bundles or an
acylindrical, hyperbolizable 3-manifold.  Let \{$N_i = {\mathbb H^3} /
\Gamma_i$\} be a sequence of hyperbolic 3-manifolds with base point in $C(N_i)$
such that each $N_i$ is homeomorphic to the interior of
$M$.  If \{$N_i$\} converges geometrically to $N = {\mathbb H^3} / \Gamma$,
then
\{$\Lambda_{\Gamma_i}$\} converges to $\Lambda_\Gamma$ in the Hausdorff
topology.
\end{corollary}

Another corollary, which does not involve base point considerations, is:

\begin{corollary} Let $M$ be a book of $I$-bundles or an acylindrical,
hyperbolizable 3-manifold.  Let $\{N_i = \HHH / \Gamma_i\}$ be a sequence
of hyperbolic 3-manifolds homeomorphic to the interior of $M$.  If
the $\{N_i\}$ converge geometrically to $N = \HHH / \Gamma$, and $\Gamma$
is nonabelian, then $\{\Lambda_{\Gamma_i}\}$ converges to
$\Lambda_\Gamma$ in the Hausdorff
topology.
\end{corollary}

In a future paper \cite{fan3}, the main result of this paper will be used to
show that for  hyperbolic 3-manifolds homeomorphic to the interior of a fixed
hyperbolizable 3-manifold with incompressible boundary, the injectivity
radius in the convex cores is uniformly bounded above.  

In \S 2 of this paper, we will introduce some background
material in hyperbolic geometry and spaces of hyperbolic 3-manifolds, and
review the relevant lemmas from
\cite{fan1}.  In \S 4, we will prove the main theorem in the case of
a book of $I$-bundles.  In \S 5, we will prove the main theorem
in the case of an acylindrical, hyperbolizable 3-manifold. In the last
section, we present some corollaries.

\section{Background Material}

\subsection{Hyperbolic Geometry and Kleinian Groups}

In this section, we review the relevant lemmas from \cite{fan1} used in the proof
of the main theorem. For the sake of brevity, we will assume that the reader is
familiar with the background material given in Section 2 of \cite{fan1}. More
details about hyperbolic geometry and Kleinian groups can be found in Beardon
\cite{beardon}, Benedetti and Petronio \cite{benedetti}, Canary-Epstein-Green
\cite{can/ep/green}, and Maskit \cite{maskit}.  

First we note that because hyperbolic space is negatively curved, injectivity radius
strictly increases while travelling out the product structure of a geometrically
finite end.

\begin{lemma} 
\label{expinj} (Lem 2.1 in \cite{fan1}) Let $N = \HHH / \Gamma$ be a hyperbolic
3-manifold. Let $U$ be a component of $N - C(N)$, and let
$S$ be the component of $\partial C_1(N)$ associated to $U$. Then
$U$ has a product structure $S \times (0,\infty)$ with nearest point
coordinates.  For $(x,t) \in S \times (0,\infty)$, the function
$inj_N(x,t)$ is strictly increasing in $t$. 
\end{lemma}

Next, we see that given an upper
bound on injectivity radius in a covering space, one can also deduce an upper bound
on injectivity radius in the base manifold.

\begin{lemma} 
\label{coverinj} (Lem 2.2 in \cite{fan1}) Let $M$ be a cover of $N$, with covering
map $p$.  Then for $x
\in M$, $inj_N(p(x)) \leq inj_M(x)$.
\end{lemma}

For a cover associated to a boundary component of the convex core of a manifold, a
component of the complement of the convex core of the base space lifts
homeomorphically to a component of the complement of the convex core of the cover.

\begin{lemma} 
\label{GFcover}  (Lem 2.3 in \cite{fan1}) Let $\delta > 0$. Let $N$ be a
hyperbolic 3-manifold.   Let $U$ be a component of
$N - int \ C(N)$, and let $S$ be
the component of
$\partial C(N)$ associated to $U$. Let $M = \HHH / \pi_1(S)$ be a cover of $N$
with projection map $p$.  Then there exists a lift $\tilde{U}$ of $U$ such
that $p|_{\tilde{U}}: \tilde{U} \rightarrow U$ is a homeomorphism, 
$(p|_{{\tilde{U}}})^{-1}(S)$ is a component of $\partial C(M)$, and
$\tilde{U}$ is a component of $M - int \ C(M)$. 
Furthermore, if $T \subset U$ is a component
of $\partial C_\delta(N)$, then $(p|_{{\tilde{U}}})^{-1}(T)$ is a
component of $\partial C_\delta(M)$.
\end{lemma}

We can also show that if a neighborhood of a geometrically finite end of a covering
space $M$ embeds as a neighborhood of a geometrically finite end in the base
manifold, then the neighborhood associated to $M - int \ C(M)$ embeds in the
base manifold as well.

\begin{lemma}
\label{GFproj}  (Lem 2.4 in \cite{fan1}) Let $0 < \epsilon < \epsilon_3$, and let
$\delta
\geq 0$.  Let
$M = \HHH / \Gamma$ be a hyperbolic 3-manifold which covers another hyperbolic
3-manifold $N = \HHH / \Gamma$ with projection map $p$.  Let
$E$ be a geometrically finite end of $M^\circ_\epsilon$, and let
$\tilde{V}$ be a component of $M - int \ C_\delta(M)$ such
that
$\tilde{V} \cap M^\circ_\epsilon$ is a neighborhood of $E$.  Let
$\tilde{S} = \partial C_\delta(M) \cap \tilde{V}$.   Suppose
there exists a neighborhood $\tilde{U}$ of $E$ such that $\tilde{U}
\subset \tilde{V}$ and $\tilde{U}$ embeds in $N$.  Then
$p(\tilde{V})$ is a component of $\NOE - C_\delta(N)$, and 
$S = p(\tilde{S})$ is a boundary component of $\partial
C_\delta(N)$. 
\end{lemma}

A nice property of topologically tame manifolds is that given a compact
set in a topologically tame manifold, it is possible to find a compact core that
contains the compact set.

\begin{lemma}
\label{compactincore} (Lem 2.6 in \cite{fan1}) Let $N$ be a topologically tame
hyperbolic 3-manifold.  Let
${\mathcal P}$ be a collection of Type I and Type II components of
$N_{thin(\epsilon)}$.  Let $K$ be a compact set in $N - {\mathcal P}$.  Then there
exists a relative compact core $R$ of
$N - {\mathcal P}$ such that $K \subset R$ and the components of $(N - \mathcal{P}) - R$ are
topologically a product.
\end{lemma}

Another property of topologically tame manifolds is that for any relative compact
core $R$ of a topologically tame manifold $N$ where the components of $\partial R -
P$ are incompressible, the components of
$(N - {\mathcal P}) - int \ R$ possess a product structure.

\begin{lemma}
\label{prodstructure} (Lem 2.7 in \cite{fan1}) Let $N$ be a topologically tame
hyperbolic 3-manifold.  Let
${\mathcal P}$ be a collection of Type I and Type II components of
$N_{thin(\epsilon)}$.  Let $R$ be a relative compact core of $N - {\mathcal P}$ with
associated parabolic locus $P$.  Let
$\{S_j\}$ be the components of
$\partial R - P$, and let $U_j$ be the component of $(N - {\mathcal P}) - int \ R$
associated to
$S_j$.  If the
\{$S_j$\} are incompressible, then the \{$U_j$\} are homeomorphic to $S_j \times [0,
\infty)$. 
\end{lemma}

Now we present a general fact about ends of hyperbolic manifolds and their covers
which will be used throughout the paper.

\begin{lemma}
\label{not in C(N)} (Lem 2.8 in \cite{fan1}) Let $N$ be a topologically tame
hyperbolic 3-manifold.  Let
$0 <
\epsilon < \epsilon_3$, and let $\delta \geq 0$.  Let ${\mathcal P}$ be a collection of
Type I components of $N_{thin(\epsilon)}$.  Let $S$ be an
incompressible separating surface of $N - {\mathcal P}$ such that $\partial S \subset
\partial {\mathcal P}$.  Let $M = \HHH / \pi_1(S)$ be a cover of $N$ with
covering map
$p$.  Let $U$ be the closure of a component of $(N - {\mathcal P}) - S$ such that
there exists
$\tilde{U} \subset M$ such that
$p|_{\tilde{U}}:\tilde{U} \rightarrow U$ is a homeomorphism and
$(p|_{{\tilde{U}}})^{-1}(S) \subset \overline{C_\delta(M)}$.  Then
$$(p|_{{\tilde{U}}})^{-1}[(N
- C_\delta(N)) \cap U] = (M - C_\delta(M)) \cap
\tilde{U}$$ 
\end{lemma}
 
\subsection{Simplicial Hyperbolic Surfaces}

In this section, we briefly review simplicial hyperbolic
surfaces.  For more details, we refer the reader to Section 3.2 of \cite{fan1}, or 
\cite{canary1}.

Let $f:
S \rightarrow N$ be a proper map from $S$ into $N$ where $S$ has triangulation
$T$. We will use the notation $f: (S,T) \rightarrow N$ to denote such a
map.  Suppose
$f$ has the following properties: 
\begin{enumerate}
\item $f$ weakly preserves parabolicity, 
\item $f$ maps every edge $e$ in $T$ to a geodesic arc, and
\item $f$ maps each face of $T$ to a non-degenerate, totally geodesic triangle
in $N$.
\end{enumerate}    Then $f$ is a {\it simplicial pre-hyperbolic surface\/}. 
Let $ang \ f(v)$ be the total angle about a vertex $f(v)$, that is, the sum of the
angles based at $f(v)$ in each of the geodesic triangles which share $f(v)$ as a
vertex.  If $ang \ f(v) \geq 2\pi$ for every internal vertex $v$, then
$f$ is a {\it simplicial hyperbolic surface\/}.  This angle condition guarantees
that the intrinsic geometry of $f(S)$ is like that of a surface of curvature
$\leq -1$.

The map $f$ induces a piecewise Riemannian metric on $S$ called the {\it
simplicial hyperbolic structure\/}, denoted $\tau$.  The surface
$(S,\tau)$ has curvature $\leq -1$ at every point except at the vertices which
have ``concentrated'' negative curvature because the total angle about each
vertex is $\geq 2\pi$.  By the Gauss-Bonnet Theorem for hyperbolic triangles,
the area of $(S,\tau)$ is
$$area(S,\tau) = 2\pi|\chi(S)|-\sum_{v\in T}(ang \ f(v)-2\pi)$$  (Lem 1.13,
Bonahon \cite{bonahon}) In particular,
$area(S,\tau) \leq 2\pi|\chi(S)|$.

Using the area bound of $(S, \tau)$, we can deduce a bound on injectivity
radius for points in the image of a $\pi_1$-injective simplicial hyperbolic surface.

\begin{lemma} 
\label{shsinj} (Lem 3.1 in \cite{fan1}) Let $f: S \rightarrow N$ be a
$\pi_1$-injective simplicial hyperbolic surface.  Then there exists a constant
$K_S$ that depends on the Euler characteristic of $S$ such that for $x \in f(S)$,
$inj_N(x) \leq K_S$.  
\end{lemma} 

Now let us explicitly describe a construction due to Bonahon \cite{bonahon} of a
simplicial pre-hyperbolic surface which will be used in the proof of the book of
$I$-bundles case.  Let
$T$ be a triangulation of
$S$ with no doubly ideal edges, and let $f: S \rightarrow N$ be a proper map that
weakly preserves parabolicity, maps every edge $e$ in $T$ with both endpoints at
the same internal vertex to a homotopically non-trivial loop in $N$, and maps no
two vertices of $T$ to the same point.  First, we homotop $f$, keeping $f(V)$
fixed, to a map $f_1: S
\rightarrow N$ such that if $e$ is a finite edge in $T$, then $f_1(e)$ is the
unique geodesic arc in its homotopy class.  Then we properly homotop $f_1$,
fixing $\bigcup f_1(e)$ over all finite edges $e$ to a map $f_2: S \rightarrow
N$ such that if $e^\prime$ is a half infinite edge, then
$f_2(e^\prime)$ is the half infinite geodesic ray which has the same endpoint
and is properly homotopic to $f(e^\prime)$.  Note that in the universal cover,
a lift of $f_2(e^\prime)$ connects a lift of $f(v)$ and a fixed point of a
parabolic element of $\pi_1(N)$.  Recall that we are assuming that there are no
doubly ideal edges, so there are no edges in the triangulation with endpoints at
two distinct ideal vertices.  Finally, we properly homotop
$f_2$, fixing $f_2(T)$, to a map $F(f, T): S \rightarrow N$ such that each face
of $T$ is taken to the totally geodesic triangle spanned by the images of its
edges.  This homotopy results in a well-defined map $F(f,T)$ which is a
simplicial pre-hyperbolic surface.  

A simplicial hyperbolic surface $f: (S,T) \rightarrow N$ is {\it useful\/} if the
triangulation $T$ contains exactly one internal vertex, and $T$ contains a
{\it distinguished edge\/} $e$ which passes through $v$ and is mapped to a closed
geodesic. It is relatively easy to check if a simplicial hyperbolic surface is
mapped into the convex core. 

\begin{lemma}
\label{shsincore} (Lem 3.7 in \cite{fan1}) Let $f:(S,T) \rightarrow N$
be a simplicial hyperbolic surface.  If $f$ maps every internal vertex of $T$
into the convex core, then $f(S) \subset C(N)$.  
In particular, if $f$ is a
useful simplicial hyperbolic surface, then $f(S) \subset C(N)$.
\end{lemma}

\subsection{Injectivity Radius Bounds for Hyperbolic $I$-bundle Convex Cores}

In this section, we state the theorems on injectivity radius bounds for hyperbolic
$I$-bundle convex cores.  The following theorem of Kerckhoff-Thurston
\cite{kerck/thur} states that given a Kleinian group that is type-preserving
isomorphic to a Fuchsian group, then there exists an upper bound on the injectivity
radius for points in the convex core of the associated hyperbolic 3-manifold.  

\begin{theorem}
\label{KT thm} (Kerckhoff-Thurston \cite{kerck/thur}) Let $\Theta$ be a cofinite
area torsion-free Fuchsian group, and let $S = \HH / \Theta$.  Then there exists a
constant $K_S$ such that for any Kleinian group $\Gamma$ such that there exists a
type-preserving isomorphism between $\Theta$ and $\Gamma$ and for $x \in C(N)$
where $N = \HHH / \Gamma$,
$inj_N(x) \leq K_S$.
\end{theorem}

A proof of their theorem also appears in Canary \cite{canary1}.  A special case of
the above theorem is as follows: \ given an
$I$-bundle
$M$ over a closed surface $S$ and a hyperbolic 3-manifold
$N$ without cusps such that $N$ is homeomorphic to the interior of $M$, then there
exists an upper bound on injectivity radius for points in the convex core of $N$. 

The following extension of their theorem will be the main tool with which we prove
our main theorem.

\begin{theorem}
\label{surface bound} (Thm 5.2  in \cite{fan1}) Let $\Theta$ be a cofinite area
torsion-free Fuchsian group, and let $S = \HH / \Theta$.  Then there exists a
constant
$L_S$ such that for any  Kleinian group $\Gamma$ such that there exists an
isomorphism between $\Theta$ and $\Gamma$ that weakly preserves parabolicity, and
for $x \in C(N)$ where $N =
\HHH / \Gamma$, $inj_N(x) \leq L_S$.
\end{theorem}

\subsection{Spaces of Hyperbolic Manifolds}

In this section, we introduce the space of
discrete, faithful representations of a Kleinian group $\Gamma$.  Given an
orientable, irreducible  3-manifold $M$, we say that
$M$ is {\it hyperbolizable\/} if there exists a discrete faithful {\it
representation\/}, $\rho: \pi_1(M) \rightarrow Isom^+(\HHH)$ such that
$N = \HHH / \rho(\pi_1(M))$ is a hyperbolic 3-manifold homeomorphic to the
interior of $M$, where $N$ has a geometric structure given by the representation. 
In this case, $N$ is a {\it hyperbolization\/} of $M$. 

A Kleinian group is {\it Fuchsian\/} if $\rho: \Gamma \rightarrow Isom^+(\HH)
\subset Isom^+(\HHH)$ so that $\Gamma$ acts properly discontinuously on $\HH
\subset \HHH$, and $\HH / \Gamma$ is a hyperbolic orbifold.  If $\Gamma$ is
torsion-free, then $\HH / \Gamma$ is a hyperbolic surface.
 
Let ${\mathcal D}(\pi_1(M))$ = {\it \{discrete faithful representations of\/}
$\pi_1(M)$ {\it into\/} $\PSLC$\}, where  ${\mathcal D}(\pi_1(M)) \subset {\rm
Hom}(\pi_1(M), \PSLC)$.  We can give ${\mathcal D}(\pi_1(M))$ the compact-open
topology, i.e.,
$\rho_i \rightarrow \rho$ if and only if $\rho_i(g) \rightarrow
\rho(g)$ for every $g \in \pi_1(M)$.  Given a representation $\rho: \pi_1(M)
\rightarrow \PSLC$,  $N_\rho = \HHH / \rho(\pi_1(M))$ is a hyperbolic 3-manifold
that is homotopy equivalent to $M$ via a homotopy equivalence induced by the map
$\rho$.  It is convenient to consider two manifolds to be equivalent if their
representations in
${\mathcal D}(\pi_1(M))$ differ by an element of $\PSLC$.  Under this
equivalence relation, we say that $N_\rho$ is a {\it marked\/} hyperbolic
manifold, and we consider $AH(M) = {\mathcal D}(\pi_1(M)) / \PSLC$ to be the space of
all marked hyperbolic 3-manifolds homotopy equivalent to
$M$.  We give $AH(M)$ the induced topology which we call the {\it algebraic
topology\/}.  If a sequence $\{\rho_i\}$ of representations converges to $\rho$
under the algebraic topology, then we say that the sequence $\{\rho_i\}$ {\it
converges algebraically\/} to the {\it algebraic limit\/} $\rho$.

We will be exclusively interested in orientable, irreducible 3-manifolds with 
incompressible boundary.  An orientable manifold is {\it irreducible\/} if every
embedded sphere bounds a ball.  An orientable manifold $M$ has {\it incompressible
boundary\/} if $i_*: \pi_1(S)
\rightarrow \pi_1(M)$ is injective for each component
$S$ of $\partial M$.   In fact, an orientable, irreducible 3-manifold $M$ has
incompressible boundary if and only if its fundamental group is {\it freely
indecomposable\/}, i.e., if $\pi_1(M) = G
\ast H$, then $G$ or $H = \{1\}$. (see Thm 7.1, Hempel \cite{hempel}) Also, $M$
has incompressible boundary if and only if $M$ contain no compressing disks.  If
$i:(D,
\partial D)
\rightarrow (M,\partial M)$ is an embedding and $i(\partial D)$ is a homotopically
non-trivial loop in $\partial M$, then we say $i(D)$ is a {\it compressing disk\/}
for $M$.  More generally, a compact surface
$S$ that is not $D^2$ or $S^2$ is {\it incompressible\/} in $M$ if $S$ is properly
embedded in
$M$ and
$i_*:
\pi_1(S)
\rightarrow \pi_1(M)$ is injective.  An orientable, irreducible 3-manifold is {\it
Haken\/} if it contains an incompressible surface.  Note that every compact,
orientable, irreducible 3-manifold with non-empty boundary is Haken. 

Furthermore, we only consider manifolds that are atoroidal.  An
embedding of a torus $i: T^2 \rightarrow M$ is {\it essential\/} if $i_*$ is
injective and $i$ is not homotopic to a map $j: T^2 \rightarrow \partial M$. 
A manifold $M$ is {\it atoroidal\/} if it contains no essential tori.  Thurston's
Geometrization Theorem states that the interior of any compact, oriented,
irreducible, atoroidal 3-manifold with nonempty boundary admits a hyperbolic
structure.

Let $M$ be a compact, orientable, irreducible 3-manifold with
incompressible boundary.  Then
$M$ is {\it acylindrical\/} if every properly embedded incompressible
annulus can be properly homotoped into the boundary of $M$.  For an acylindrical
manifold, Thurston \cite{thurston2} has shown that its representation space is
compact:

\begin{theorem}
\label{AH(M)cmpt} (Thurston's Compactness Thm \cite{thurston2})  Let $M$ be
an  atoroidal, acylindrical 3-manifold.  Then $AH(M)$ is compact. 
\end{theorem}

We can also consider the space of Kleinian groups with the {\it geometric
topology\/} or {\it Chabauty topology\/}, that is, the topology of closed
subgroups.  A sequence $\{\Gamma_i\}$ of Kleinian groups {\it converges
geometrically\/} to a group $\Gamma$ if the following conditions are satisfied:

\begin{enumerate}
\item if $g \in \PSLC$ is an accumulation point of a sequence $\{g_{i}\}
\in \{\Gamma_i\}$, then $g \in \Gamma$, and
\item if $g \in \Gamma$, then there exists a sequence $\{g_i\}$ such that each $g_i
\in \Gamma_i$ and $g_i \rightarrow g$.
\end{enumerate} 

One can also think of geometric convergence in terms of hyperbolic
3-manifolds with base frame.  Choose a base point $p$ in $\HHH$ and an orthonormal
frame $\omega$ in the tangent space at $p$.  Given a manifold
with base frame $(M,e)$, there exists a unique Kleinian group $\Gamma$ with the
property that
$\Gamma$ is the group of covering transformations acting on $\HHH$ such that
$(\HHH / \Gamma, \omega) = (M,e)$, where $e$ is the image of the standard frame
$\omega$ at $p$.  

We can now give the space of hyperbolic 3-manifolds with base frame the Chabauty
or geometric topology.  We say $f$ is a {\it framed\/} $(K,r)$-{\it approximate
isometry\/} between two manifolds with base frame $(M_1,e_1)$  and $(M_2,e_2)$
if $f: (X_1,e_1) \rightarrow (X_2,e_2)$ is a diffeomorphism such that
$B_{M_1}(x_1,r) \subseteq (X_1,x_1) \subseteq (M_1,x_1)$ , $B_{M_2}(x_2,r)
\subseteq (X_2,x_2) \subseteq (M_2,x_2)$, $Df(e_1) = e_2$, and
$$\frac{d(x,y)}{K} \leq d(f(x),f(y)) \leq K d(x,y)$$ for all $x,y \in X_1$. So
$\{(M_i,e_i)\}$ {\it converges geometrically\/} to $(M,e)$ if there exists
a sequence of $(K_i,r_i)$-approximate isometries \{$f_i:(M_i,e_i) \rightarrow
(M,e)$\} such that
$K_i \rightarrow 1$ and $r_i \rightarrow \infty$. 
The topology induced by framed $(K,r)$-approximate isometries coincides with the
geometric topology.  (Cor 3.2.11, Canary-Epstein-Green \cite{can/ep/green}) 
Thus, a sequence of manifolds with base frame $(M_i,e_i)$ converges to $(M,e)$
in the geometric topology if and only if their corresponding Kleinian groups
$\Gamma_i$ converge to $\Gamma$ in the geometric topology.

\section{The Book of I-Bundles Case}

In this section, we will prove the main theorem in the case that $M$ is a book
of $I$-bundles. First we will give a sketch of the proof of the theorem for a
simpler example of a book of $I$-bundles.  Then we will prove the main theorem
for a general book of $I$-bundles. 

\subsection{The Motivating Example}

Let us consider a simpler example of a
book of $I$-bundles called a {\it plain book of I-bundles\/}.  To construct a
plain book of $I$-bundles, let \{$B_i$\} be a collection of surfaces,
each of which is a closed orientable surface minus an open disk.  For each
$i$, let $E_i = B_i \times [0,1]$, and let
$\partial_\circ E_i$ be the annulus $\partial B_i \times [0,1]$.  Consider a
solid torus $V = D^2
\times S^1$ whose boundary is decomposed into a union of disjoint parallel
annuli, $A_1, A_1^\prime, A_2, A_2^\prime, \ldots, A_n, A_n^\prime$, such that
each annulus is homotopy equivalent to a $(1,0)$-curve on the boundary of the
torus.  We can order the annuli {\it mod
n\/}.  Let $M$ be the union of $V$ and the \{$E_i$\}, where each
$\partial_\circ E_i$ is glued to $A_i$ by a homeomorphism.
We can think of $V$ as
the binding and the \{$E_i$\} as the pages in this book of
$I$-bundles.

\begin{theorem}
\label{plain case} Let $M$ be a plain book of $I$-bundles. Then there exists a
constant $L$ such that if $N$ is a hyperbolic 3-manifold homeomorphic
to the interior of $M$ and if $x \in C(N)$, then $inj_N(x) \leq L$, where $L$
depends on the maximum genus of the boundary components of $M$.
\end{theorem}

\begin{proof} For the sake of simplicity, let us assume that $N$ is convex
co-compact.  Let us start with an outline of the sketch of the proof in this case. 
The idea is to choose a 2-complex $D \subset int \ M$ which will be a union of surfaces
\{$T_i$\} which are isotopic to the boundary components \{$S_i$\} of $M$. We
will then consider a map
$f:D \rightarrow N$ such that $f(D) \subset C(N)$ and such that $f$ is in the
appropriate homotopy class so that each of the components of $C(N) -
f(D)$ will lie in the image of the convex core of a manifold whose
convex core has bounded injectivity radius, thus proving the theorem.

Now we will give a more complete sketch of the proof.  First we construct a 2-complex
$D$ in $int \ M$.  For each $i$, let ${\mathcal A}_i
\subset V$ be an annulus whose boundary components are the core curve of $V$ and the
core curve of $A_i$.  Let $D$ be the 2-complex formed by gluing the boundaries
of $\{(B_i, \frac{1}{2}) \subset E_i\}$ to those of $\{{\mathcal A}_i\}$.  Then by
construction, the inclusion of $D$ into $M$ is a homotopy equivalence.  Let $S_i =
(B_i,0) \cup A_i^\prime \cup (B_{i+1},0)$ be the $i$-th boundary component of $M$, and
let
$T_i = (B_i,\frac{1}{2})
\cup {\mathcal A}_i \cup {\mathcal A}_{i+1} \cup (B_{i+1},\frac{1}{2})$, glued along
their boundaries.  Then $S_i$ is isotopic to $T_i$, and the component of $M -
D$ containing $S_i$ is homeomorphic to $S_i \times [0,1]$.

Let $M$ be our plain book of $I$-bundles described above.  Let $\rho \in {\mathcal
D}(\pi_1(M))$, and let $N = \HHH / \rho(\pi_1(M))$ be a hyperbolic 3-manifold
homeomorphic to $int \ M$, via a map $h: int \ M \rightarrow N$. 
 Let $M_i = {\mathbb H}^3 / \rho(\pi_1(S_i))$ be
a cover of $N$ with covering map
$p_i$.  

Let the vertex set $W$ on $D$ consist of exactly one internal vertex on the core
curve of $V$.  Triangulate $(D,W)$.  Let $f: D
\rightarrow N$ be a map such that for each
$i$, $f|_{T_i}: T_i \rightarrow N$ is a useful simplicial hyperbolic surface so that
$f(D) \subset C(N)$.  Moreover, construct $f$ so that it is homotopic to $h|_D$ and
hence $\pi_1$-injective.  We further require that for each $i$, 
the map $f|_{T_i}$ lifts to $\tilde{f}_i: T_i \rightarrow M_i$ such that
$\tilde{f}_i(T_i) \subset C(M_i)$.  

By Lemma \ref{GFcover}, a boundary component $A$ of $C(N)$ lifts to a boundary
component $\tilde{A}$ of $C(M_i)$. We can construct a homotopy $H_i: S_i \times [0,1]
\rightarrow M_i$ between $\tilde{A}$ and
$\tilde{f}_i(T_i)$ such that $H_i(S_i,[0,1]) \subset C(M_i)$.  Because the image of
the homotopy lies in the convex core of a surface group, by Theorem
\ref{surface bound}, for $x \in H_i(S_i,[0,1])$, $inj_{M_i}(x) \leq L_{S_i}$ where
$L_{S_i}$ depends on $S_i$.  So by Lemma \ref{coverinj}, for $x \in \bigcup
p_i(H_i(S_i,[0,1]))$,
$inj_N(x) \leq L_{S_i}$.  With some work, we can show that $$C(N) \subset
\bigcup p_i(H_i(S_i,[0,1])) \subset \bigcup p_i(C(M_i))$$  Thus, for $x \in
C(N)$, $inj_N(x) \leq max\{L_{S_i}\}$. 
\end{proof}

\subsection{The General Case}

Now we will prove the main theorem in the case of a general book of
$I$-bundles.  A {\it book of I-bundles\/} is a compact, connected,
irreducible 3-manifold with boundary $M = E \cup V$ such that
\begin{enumerate}
\item $E$ is an $I$-bundle over $B$, a non-empty compact 2-manifold with
boundary,
\item each component of $V$ is homeomorphic to $D^2 \times S^1$, 
%
%$T^2 \times [-1,1]$, CAN YOU EVEN HANDLE THIS CASE?
%
\item the set $A = E \cap V$ is the inverse image of $\partial B$ under the
bundle projection $b: E \rightarrow B$, and
\item each component of $A$ is an annulus in $\partial V$ which is homotopically
non-trivial in $V$.
\end{enumerate}
Equivalently, a compact, connected, irreducible 3-manifold $M$ is a {\it
book of $I$-bundles\/} if there exists a disjoint collection $A$ of incompressible
annuli such that each component of the manifold obtained by cutting $M$ along
$A$ is either a solid torus, or an $I$-bundle $R$ over a surface of negative
Euler characteristic such that $\partial R \cap \partial M$ is the associated
$\partial I$-bundle. 

The proof of the motivating example was simpler for several reasons: \ $N$ was convex
co-compact; the annuli in $\partial V$ were homotopic to $(1,0)$-curves, not
$(p,q)$-curves; $M$ did not contain multiple bindings, i.e., $V$ had only one
component; and there were no twisted $I$-bundles.  These are all considerations that
we will have to take into account when we prove the main theorem in the case of a
general book of
$I$-bundles.

\begin{theorem}
\label{bookresult}  Let $M$ be a book of $I$-bundles. Then there exists a
constant $L$ such that if $N$ is a hyperbolic 3-manifold homeomorphic
to the interior of $M$ and if $x \in C(N)$, then $inj_N(x) \leq L$, where $L$
depends on the maximum genus of the boundary components of $M$.
\end{theorem}

\begin{proof}  The outline of the proof is
the similar to that for a plain book of $I$-bundles.  The idea is to use the image of
a 2-complex $D$ to divide $C(N)$ into portions, each of which lies in the image of the
convex core of a manifold whose convex core has bounded injectivity radius. 

In the following lemma we will choose a 2-complex $D \subset int \ M$ such that the
inclusion of $D$ into $M$ is a homotopy equivalence, and we will construct maps
\{$g_i: S_i
\times [0,1] \rightarrow M$\} whose images can
be glued along $D$ to form a 3-manifold homeomorphic to $M$.

\begin{lemma}
\label{what is D} Let \{$S_i$\} be the boundary components of $M$.  Then there
exists a 2-complex $D \subset int \ M$ such that $D$ is a deformation retract of $M$. 
Furthermore, there exist maps \{$g_i: S_i
\times [0,1] \rightarrow M$\} with the following properties: 
\begin{enumerate}
\item for each $i$, the image $g_i(S_i,0) = S_i$,
\item the map $g_i|_{S_i \times [0,1)}: S_i \times [0,1) \rightarrow M$ is a
homeomorphism onto its image, 
\item the image $g_i(S_i,1) \subset D$, and
\item the manifold $M$ is homeomorphic to the quotient space $\coprod
(S_i,[0,1])/\sim$, where
$(x_i,1) \sim (y_j,1)$ if and only if $(g_i(x_i),1) = (g_j(y_j),1)$.
\end{enumerate}
\end{lemma}

\begin{proof}  Recall that $M = E \cup V$, where $E$ is an $I$-bundle over $B$, and
$B$ is a disjoint collection
\{$B_s$\} of surfaces with boundary.   Here, $V$ is a disjoint collection $V_1,
\ldots, V_m$ of solid tori.  The boundary of each $V_j$ consists of a collection of $2
n_j$ parallel annuli $A_{j1}, A_{j1}^\prime, A_{j2}, A_{j2}^\prime, \ldots,
A_{jn_j}, A_{jn_j}^\prime$, which we can order cyclically
$mod \ n_j$.  Let $b: E \rightarrow B$ be the bundle projection map.  Given a component
$B_s$ of
$B$, for each component
$\hat{Q}$ of $\partial B_{s}$, $b^{-1}(\hat{Q})$ is identified with some annulus
$A_{jk}$ on the boundary of some $V_j$.  In fact, there is a one-to-one
correspondence between the set \{$A_{jk}$\} and the set \{$b^{-1}(\hat{Q}) : \
\hat{Q}$ is a boundary component of some $B_s$\}. 

We define two annuli $A_{jk}$ and $A_{jl}$ to be {\it adjacent on
$\partial V_j$\/} if they are ordered consecutively on $\partial V_j$. Note that if
$A_{jk}$ and $A_{j(k+1)}$ are adjacent annuli on $\partial V_j$, then they will be
physically separated by the parallel annulus $A_{jk}^\prime$ to which no
$I$-bundle is identified.
Given an $I$-bundle $E_s$ over $B_s$, either $E_{s} = B_{s} \times [0,1]$, or $E_s$
is a twisted $I$-bundle over $B_s$.  In the latter case, $E_s = \tilde{B}_s \times
[0,1] / (x,t) \sim (\tau(x),1-t)$, where $\tau: \tilde{B}_s \rightarrow \tilde{B}_s$ is
a free involution and
$\tilde{B}_{s} / \tau = B_{s}$.  For a component
$B_s$ of
$B$, let
$\bar{B}_{s}$ denote the {\it middle surface\/} of the $I$-bundle over $B_s$, where
$\bar{B}_{s} = (B_{s},
\frac{1}{2})$ if
$E_{s} = B_{s} \times [0,1]$, and $\bar{B}_{s} = (\tilde{B}_{s}, \frac{1}{2}) /
\tau$ if $E_{s}$ is a twisted $I$-bundle over $B_{s}$. Let
$\bar{B}$ denote the set of all the middle surfaces \{$\bar{B}_s$\}.  

Let $\bar{M} = E \cup (\partial V \times [0,1])$ be a compact 3-manifold, where $E$
is an $I$-bundle over $B$ and $V$ is a disjoint collection of
solid tori as before.  For each $j$, $(\partial V_j,0)$ is
decomposed into parallel annuli $(A_{j1},0), (A_{j1}^\prime,0), 
\ldots, (A_{jn_j},0), (A_{jn_j}^\prime,0)$.  Given a component 
$B_s$ of $B$ with boundary component $\hat{Q}$, if $b^{-1}(\hat{Q})$ is identified
with $A_{jk}$ on $\partial V_j$ in $M$, then identify $b^{-1}(\hat{Q})$ with
$(A_{jk},0)$ on $(\partial V_j,0)$ in $\bar{M}$.  Consider the quotient map $q: \bar{M}
\rightarrow M$ defined such that $q|_B$ is the identity map, $q|_{\partial V \times
[0,1)}$ is a homeomorphism, and $q$ identifies $(\partial V_j,1)$ to the core curve
$\gamma_j$ of
$V_j$.   

Now we will construct a disjoint union of surfaces, $\bar{D}$, such that
$q(\bar{D})$ will be the desired 2-complex $D$. For a fixed annulus $(A_{jk},0)$ in
$(\partial V_j,0)$, construct the annulus $\bar{\mathcal A}_{jk} = (\partial \bar{B}_{s}
\cap A_{jk}) \times [0,1] \subset (\partial V_j \times [0,1])$.  Then for fixed $s$,
let
$\bar{\mathcal A}_s = \bigcup \bar{\mathcal A}_{jk}$ where the union is taken over all
components of
$(\partial \bar{B}_{s} \cap \partial V,0)$.  Then $\bar{F}_s = \bar{B}_{s} \cup
\bar{\mathcal A}_s$ is a surface with boundary \{$(\partial \bar{B}_s \cap
A_{jk},1)$\}.  Then we define
$\bar{D} = \bigcup \bar{F}_s$, and $D = q(\bar{D}) \subset int \ M$.  

%Let $h_{js}:
%\hat{B}_{js} \times [0,\frac{1}{2}]
%\rightarrow \bar{B}_{js}$ be the natural bundle projection onto the middle surface
%where $h_{js}(\hat{B}_{js},0) = \hat{B}_{js}$ and
%$h_{js}(\hat{B}_{js},\frac{1}{2}) =
%\bar{B}_{js}$.  Here
%$\hat{B}_{js}$ is the {\it top surface\/} and is defined
%$\hat{B}_{js} = B_{js}$ if
%$E_{js} = B_{js}
%\times [0,1]$, and $\hat{B}_{js} = \tilde{B}_{js}$ if $E_{js}$ is a twisted
%$I$-bundle over
%$B_{js} = \tilde{B}_{js} / \tau$.  

There is a natural product structure on each component of
$E - \bar{B}$.  For adjacent annuli $(A_{jk},0)$ and $(A_{j(k+1)},0)$, we can extend
this product structure so that the component of
$(\partial V_j \times [0,1]) - [\bar{\mathcal A}_{jk} \cup \bar{\mathcal A}_{j(k+1)} \cup
(\partial V_j,1)]$ containing
$(A_{jk}^\prime,0)$ also has a product structure as shown in Figure
~\ref{fig:prodstructure}.  This product structure induces a product structure on all
of $M - D$.  Note that as a result, we can deformation retract $M$ along the product
structure of
$M - D$ onto
$D$.

\begin{figure}[ht]
\begin{center}
\BoxedEPSF{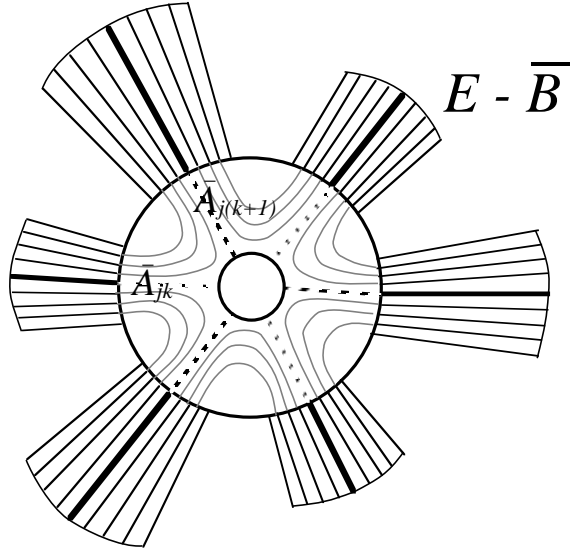 scaled 600}
\caption{Cross Section of Product Structure on $\partial V_j \times [0,1]$ in
$\bar{M}$}
\label{fig:prodstructure}
\end{center}
\end{figure}

Now let us construct the maps \{$g_i: S_i \times [0,1] \rightarrow M$\}.  Consider a
boundary component $S_i$ of $M$.  Then for each $i$, let
$g_i: S_i \times [0,1) \rightarrow M - D$ be a homeomorphism along the product
structure of $M - D$ where $g_i(S_i,0) = S_i$.  We can let $g_i(x_i,1) = lim_{t
\rightarrow 1} g_i(x_i,t)$ for $x_i \in S_i$.  Then for all $i$, $g_i(S_i,1)  \subset
D$.  By identifying the images of the maps in $D$, we see that $M$ is homeomorphic to
the quotient space
$\coprod (S_i,[0,1])/\sim$, where
$(x_i,1) \sim (y_j,1)$ if and only if $(g_i(x_i),1) = (g_j(y_j),1)$.   \end{proof}

Let $N$ be a hyperbolic 3-manifold homeomorphic to $int \ M$ via a homeomorphism $h:
int \ M \rightarrow N$.  For ease of exposition, let $\hat{g}_i: S_i \rightarrow D
\subset M$ be defined by $\hat{g}_i(x_i) = g_i(x_i,1)$ for $x_i \in S_i$.  In the next
lemma, we will construct a map $f: D \rightarrow N$ such that for each $i$, the image
$f \circ \hat{g}_i(S_i)$ lifts to the $\epsilon$-neighborhood of the convex core of a
cover of $N$.

\begin{lemma}
\label{map D in}  Let $N$ be a hyperbolic 3-manifold homeomorphic to $int \ M$ via a homeomorphism $h:
int \ M \rightarrow N$.  Let $i: D \rightarrow int \ M$ be the inclusion map, and let
$k = h
\circ i: D \rightarrow N$.  Let $\Gamma_i = (k \circ \hat{g}_i)_*[\pi_1(S_i)]$, and
let
$M_i =
\HHH / \Gamma_i$.  Then there exists a map $f: D \rightarrow N$ such that for each $i$,
$f
\circ \hat{g}_i$ lifts to a map $\widetilde{f \circ \hat{g}_i}: S_i
\rightarrow M_i$ such that $\widetilde{f \circ \hat{g}_i}(S_i) \subset
C_{\epsilon}(M_i)$.
\end{lemma}

\begin{proof} Recall that $D = q(\bar{D})$ where
$\bar{D} = \bigcup \bar{F}_s$, a disjoint union of surfaces with boundary.  We will
define the map
$f: D
\rightarrow N$ by constructing a map $\bar{f}: \bar{D}
\rightarrow N$, where for each $\bar{F}_s \subset \bar{D}$, the map
$\bar{f}|_{\bar{F}_s}$ is a simplicial hyperbolic surface that factors through the
quotient map $q$.

Let us construct a triangulation on $\bar{D}$ that will induce a ``triangulation''
on $D$.  We put ``triangulation'' in quotes, because $D$ is a 2-complex, not a surface.
First choose a vertex set
$W$ on
$D$ to be exactly one internal vertex
$v_j$ on each $\gamma_j$.   Let
$\bar{W} = q^{-1}(W)$ be the vertex set on $\bar{D}$.  Include the arcs
\{$(\partial\bar{B}_s \cap A_{jk},1) - \bar{W}$\} in the edge set of the triangulation
on $(\bar{D}, \bar{W})$.  For a fixed $s$ and
$j$, if $(\partial \bar{B}_s \cap A_{jk},0)$ is homotopic to a $(p_j,q_j)$-curve on
$(\partial V_j,0)$, then there are $p_j$ vertices on $(\partial \bar{B}_s \cap
A_{jk},1)$, and the edge set
\{$(\partial \bar{B}_s \cap A_{jk},1) - \bar{W}$\} in $\bar{D}$ is a $p_j$-to-1
cover of $(\gamma_j - W)$ in $D$.  Triangulate the remainder of $(\bar{D}, \bar{W})$,
and call it
$\bar{T}$.  Then $T = q(\bar{T})$ is a ``triangulation'' on  $(D,W)$.  

Because $D$ is a deformation retract of
$M$, the inclusion map $i: D
\rightarrow int \ M$ is a homotopy equivalence.  Then $k = h \circ i: D \rightarrow N$
is
$\pi_1$-injective.  Now we will construct a map $\hat{k}: D \rightarrow N$ such that
$\hat{k}$ will be homotopic to $k$, and $\bar{f}$ will be a
simplicial hyperbolic surface which is a ``straightening'' of
$\hat{k} \circ q$. 

For each $j$, $k(\gamma_j)$ represents either a hyperbolic element or a
parabolic element in
$\pi_1(N)$.  If $k(\gamma_j)$ represents a hyperbolic element, then its geodesic
representative,
$k(\gamma_j)^*$, lies in $C(N)$.  Let $\hat{k}: D \rightarrow N$ be a map that is
homotopic to $k: D \rightarrow N$ such that 
$\hat{k}(v_j)$ is a point on $k(\gamma_j)^*$.

If $k(\gamma_j)$ represents a parabolic element, then $k(\gamma_j)$ has no geodesic
representative in $N$.  Normalize so that the fixed point of the parabolic element
which $k(\gamma_j)$ represents is at infinity in the upper half space model of
$\HHH$.  Recall that $\gamma_j$ is the core curve of $V_j$.  Consider
the collection
$\{S_i\}$ of boundary components of
$M$ such that
$S_i
\cap V_j \neq \emptyset$.  Consider the collection of covers \{$M_i = \HHH /
\Gamma_i$\} of $N$ associated to the  $\{S_i\}$,  where each $\Gamma_i = (k \circ
\hat{g}_i)_*[\pi_1(S_i)]$. Because each $C(M_i)$ is non-empty, the convex hull of
each
$\Lambda_{\Gamma_i}$ must also contain a geodesic ray
$Y_i$ with endpoint at infinity.  Consider a
horoball $L_j = \{(z,t): t \geq c_j\}$ about infinity.  Under the action of the
normalized parabolic element $z
\mapsto z + 1$, the horoball $L_j$ induces a rank one cusp in
$N$, the boundary of which is an infinite annulus.  Then for large enough $c_j$, there
exists an arc
$\tilde{\alpha}_j$ on $\partial L_j$ such that
$\alpha_j =
\tilde{\alpha}_j / (z
\mapsto z + 1)$ is a closed curve in $N$ on the infinite annulus $\partial L_j / (z
\mapsto z + 1)$, and such that $\tilde{\alpha}_j$ is contained in the
$\epsilon$-neighborhood of each $Y_i$ so that $\tilde{\alpha}_j \subset
CH_\epsilon(\Lambda_{\Gamma_i})$ for each
$\Gamma_i$.  Let $\hat{k}: D \rightarrow N$ be a map that is homotopic to $k: D
\rightarrow N$ such that $\hat{k}(v_j)$ is a point on $\alpha_j$.
 
Using Bonahon's construction, let $\bar{f} = F(\hat{k} \circ q,T): \bar{D}
\rightarrow N$.  Then for each $\bar{F}_s \subset \bar{D}$, $\bar{f}|_{\bar{F}_s}$ is a
simplicial pre-hyperbolic surface.  Because
$\bar{f}$ is in the homotopy class of $\hat{k} \circ q$,  $\bar{f}$ is a
$\pi_1$-injective map.  By construction, the map $\bar{f}$ respects the quotient map
$q$, so there exists a map $f: D \rightarrow N$ which is
a simplicial pre-hyperbolic 2-complex, i.e., $f$ weakly preserves
parabolicity, maps every edge in
$T$ to a geodesic arc, and maps each face of $T$ to a
non-degenerate totally geodesic triangle in $N$.  Furthermore, by construction, $f$ is
homotopic to
$k: D
\rightarrow N$ and hence is $\pi_1$-injective. 

Now we will show that $f \circ \hat{g}_i(S_i)$ lifts to $C_\epsilon(M_i)$.  Because
$f \circ \hat{g}_i(S_i)$ is homotopic to $k \circ \hat{g}_i(S_i)$, by
the Lifting Theorem,  $f \circ \hat{g}_i(S_i)$ lifts to a map $\widetilde{f \circ
\hat{g}_i}: S_i \rightarrow M_i$.  Because $\hat{g}_i(S_i)$ is a subset of $D$, $f
\circ
\hat{g}_i$ is a simplicial pre-hyperbolic surface.  Then
$\widetilde{f \circ \hat{g}_i}: S_i \rightarrow M_i$ is also a simplicial
pre-hyperbolic surface.  Because $C_\epsilon(M_i)$ is also convex, by an argument
similar to that in the proof of Lemma
\ref{shsincore} \cite{fan1}, in order to show that $\widetilde{f \circ \hat{g}_i}(S_i)
\subset C_\epsilon(M_i)$, it suffices to show that for each internal vertex $v_j$ in
the triangulation on
$\hat{g}_i(S_i)$, $f(v_j)$ lifts to $C_{\epsilon}(M_i)$.  

Let $v_j \in \gamma_j$ be an internal vertex in the triangulation on
$\hat{g}_i(S_i) \subset D$.  For each $\gamma_j$, $k(\gamma_j)$ represents either a
hyperbolic element or a parabolic element in $\Gamma_i$.   If
$k(\gamma_j)$ represents a hyperbolic element in $\Gamma_i$, then $f(v_j) =
\hat{k}(v_j) \in k(\gamma_j)^*$.  The closed geodesic $k(\gamma_j)^*$ lifts to
$\HHH$ to an axis of a hyperbolic element of $\Gamma_i$ which lies in
$CH(\Lambda_{\Gamma_i})$.  Then
$f(v_j)$ lifts to
$C(M_i)$.  If $k(\gamma_j)$ represents a parabolic element in $\Gamma_i$, then, by
construction,
$f(v_j)$ lifts to $\HHH$ to a point on $\tilde{\alpha}_j$ which lies in
$CH_\epsilon(\Lambda_{\Gamma_i})$.  Therefore, $f(v_j)$ lifts to
$C_\epsilon(M_i)$, and  we have shown that $\widetilde{f \circ \hat{g}_i}(S_i) \subset
C_\epsilon(M_i)$.
\end{proof}

The next step in the proof of the book of $I$-bundles case involves finding a compact
core
$R$ of
$N$, such that the points in
$R$ have bounded injectivity radius.

Given $0 < \epsilon < \epsilon_3$. By the results of McCullough
\cite{mccullough} and Kulkarni-Shalen
\cite{kul/shalen} and Lemma \ref{compactincore}, there exists a compact core $R^*$ of
$\overline{C_\epsilon(N)}$ such that  $\partial C_\epsilon(N) \cap
C^\circ_\epsilon(N)
\subset \partial R^*$, and such that $f(D) \subset R^*$.
 
In the next lemma, we will choose a compact core $R$ of $N$ such that
$R$ contains $R^*$, $R$ lies in $C_{2 \epsilon}(N)$, and such that points in $R$
have bounded injectivity radius.

\begin{lemma}
\label{R inj bdd} Let $\Gamma_i = (k \circ
\hat{g}_i)_*[\pi_1(S_i)]$, and let $M_i = \HHH / \Gamma_i$ be a cover of $N$ with
covering map $p_i$.  Let \{$U_i$\} be the
components of $C_{2 \epsilon}(N)
 - int \ R^*$.   For all $0 < \epsilon < \epsilon_3$, there exists a compact core
$R$ of $N$ bounded by the surfaces $\{T_i\}$ such that the
following are true: 

\begin{enumerate}
\item $R \subset C_{2 \epsilon}(N)$, 
\item if $x \in R$, then $inj_N(x) \leq max \{L_{S_i} + 2 \epsilon\}$, and
\item for each $i$, there exists $\tilde{U}_i \subset M_i$ such that
$p_i|_{\tilde{U}_i}: \tilde{U}_i \rightarrow U_i$ is a homeomorphism, and
$(p_i|_{\tilde{U}_i})^{-1}(T_i) \subset \overline{C_{2 \epsilon}(M_i)}$, 
\end{enumerate}
\end{lemma}

\begin{proof} Let us start with a sketch of the
proof.  First we will construct surfaces \{${\mathcal T}_i \subset C_{2 \epsilon}(N)$\}
such that for each $i$, there exists a homotopy between $f \circ \hat{g}_i:S_i
\rightarrow N$ and the inclusion map $i|_{{\mathcal T}_i}:{\mathcal T}_i \rightarrow N$.  We
will construct each homotopy so that its image lies in the image of the convex core of
a manifold whose convex core has bounded injectivity radius.  Then we will construct
surfaces \{$T_i \subset C_{2 \epsilon}(N)$\} such that there exists a homotopy
between the inclusion maps of ${\mathcal T}_i$ and $T_i$.  Again, we will
construct each homotopy so that its image lies in the image of the convex core of a
manifold whose convex core has bounded injectivity radius.  Finally, we will show that
the \{$T_i$\}  bound a compact core $R$ and that $R$ is contained in the union of the
images of these homotopies.  Thus, by construction, all points in $R$ will have bounded
injectivity radius. 

Note that $\pi_1(N) = \pi_1(\overline{C_{\epsilon}(N)})$, so that $R^*$ is also a
compact core of $N$.  Because $N$ is homeomorphic to
$int \ M$, and $R^*$ is a compact core of $N$, we can conclude that $M$ is homeomorphic
to $R^*$.  (Thm 1, McCullough-Miller-Swarup
\cite{mmswarup})  Then because $M$ has incompressible boundary, 
$R^*$ also has incompressible boundary.  Note that $R^*$ is also a
compact core for
$C_{2 \epsilon}(N)$ where the components of $\partial R^*$ are incompressible in $C_{2
\epsilon}(N)$.   Let $\{Q_i\}$ be the components of $\partial R^*$.   Then by Lemma
\ref{prodstructure}, each component $U_i$ of $C_{2 \epsilon}(N)
 - int \ R^*$ possesses a product
structure $Q_i \times [0, \infty)$.  Because $i_*(\pi_1(Q_i)) = \Gamma_i$,
by the Lifting Theorem, the inclusion map $i: Q_i \times [0,\infty)
\rightarrow U_i$ lifts to a map $\tilde{i}: Q_i \times [0,\infty) \rightarrow
M_i$.  Let $\tilde{U}_i = \tilde{i}(Q_i,[0,\infty))$.  Then $\tilde{U}_i$
is also possesses a product structure $Q_i \times [0, \infty)$, and
$p_i|_{\tilde{U}_i}:\tilde{U}_i \rightarrow U_i$ is a homeomorphism.

Now we will show that for each $i$, $\tilde{U}_i \subset C_{2 \epsilon}(M_i)$.  Suppose
not.  Then there exists $\tilde{x} \in \tilde{U}_i - C_{2 \epsilon}(M_i)$ and a
component
$\tilde{G}_i$ of $\partial C_{2 \epsilon}(M_i)$ such that $\tilde{G}_i$
separates $\tilde{x}$ from $C_{2 \epsilon}(M_i)$.  Consider a geodesic ray
$\bar{g}_{\tilde{x}}$ in $M_i - C_{2 \epsilon}(M_i)$ that is perpendicular to
$\tilde{G}_i$ and passes through $\tilde{x}$.  Let
$g_{\tilde{x}}$ be the portion of $\bar{g}_{\tilde{x}}$ beginning at $\tilde{x}$. 
There exists a constant $\eta > 0$ such that $\tilde{x} \in (M_i)^\circ_\eta$.  Recall
that by Lemma
\ref{expinj}, the injectivity radius strictly increases out a geometrically finite
end, therefore the ray $g_{\tilde{x}}$ is contained in
$(M_i)^\circ_\eta$.  Without loss of generality, let us assume that
$g_{\tilde{x}}$ intersects $(Q_i,0)$ transversely.  Then one of the following three
cases occurs: \ $g_{\tilde{x}}$ intersects
$(Q_i,0)$ an odd number of times,
$g_{\tilde{x}}$ intersects a different boundary component of the closure of
$\tilde{U}_i$, or all but a compact portion of 
$g_{\tilde{x}}$ is contained in $\tilde{U}_i$.

Suppose $g_{\tilde{x}}$ intersects $(Q_i,0) \subset \tilde{U}_i$ an odd number of
times.  Because
$(Q_i,0) \subset \tilde{U}_i$ is an embedded surface separating
$M_i$, $(Q_i,0)$ has two sides.  Let the component of
$M_i - (Q_i,0)$ containing
$\tilde{U}_i$ be the positive side of
$(Q_i,0)$.  Then the last time $g_{\tilde{x}}$ intersects $(Q_i,0) \subset
\tilde{U}_i$,
$g_{\tilde{x}}$ passes from the positive to the negative side of 
$(Q_i,0)$.  

We can also let the side of $C_{2 \epsilon}(N) - (Q_i,0)$ containing $U_i$ be the
positive side of $(Q_i,0) \subset U_i$.  Recall that because
$f \circ g_i(S_i,1) \subset R^*$ and $U_i$ is a component of $C_{2
\epsilon}(N) - int \ R^*$, $f \circ g_i(S_i,1)$ lies on the negative
side of
$(Q_i,0) \subset U_i$.  

Because $p_i|_{\tilde{U}_i}$ is a homeomorphism and $(Q_i,0) \subset U_i$ lifts to
$(Q_i,0) \subset \tilde{U}_i$, $\widetilde{f \circ g_i}(S_i,1)$ lies to the negative
side of $(Q_i,0) \subset \tilde{U}_i$.  Note that $\widetilde{f \circ g_i}(S_i,1)$
is an incompressible separating surface of $M_i$, and that
$\widetilde{f \circ g_i}:(S_i,1) \rightarrow M_i$ is homotopic to the
inclusion map $i: (Q_i,0) \rightarrow M_i$.  Because
$g_{\tilde{x}}$ leaves every compact set of $M_i$ and all but a compact portion of
$g_{\tilde{x}}$ lies in $M_i - \tilde{U}_i$, if
$g_{\tilde{x}}$ intersects $(Q_i,0)
\subset \tilde{U}_i$, then $g_{\tilde{x}}$ must also intersect
$\widetilde{f \circ g_i}(S_i,1)$.  But $g_x \subset M_i - C_{2 \epsilon}(M_i)$ and
$\widetilde{f \circ g_i}(S_i,1) \subset C_{\epsilon}(M_i)$.  So this is a
contradiction.  

Suppose $g_{\tilde{x}}$ intersects a different boundary component of the closure of
$\tilde{U}_i$.  Then $p_i(g_{\tilde{x}})$ intersects some
component $H_i$ of $\partial C_{2 \epsilon}(N)$ such that $H_i$ is a boundary
component of the closure of $U_i$ in $N$.  By Lemma \ref{GFcover},
$H_i$ lifts to a component $\tilde{H}_i$ of $\partial C_{2 \epsilon}(M_i)$.  Because
$g_{\tilde{x}}$ intersects $\tilde{H}_i$, there is a portion of 
$\bar{g}_{\tilde{x}}$ that contains $\tilde{x}$ and joins two components of $\partial
C_{2 \epsilon}(M_i)$.  Because $C_{2
\epsilon}(M_i)$ is convex, $\tilde{x} \in g_{\tilde{x}}$ lies in $C_{2
\epsilon}(M_i)$, which is a contradiction.

Then all but a compact portion of $g_{\tilde{x}}$ is contained in
$\tilde{U}_i$.  Let $\tilde{W}_i$ be the closure of the component of
$M_i - \tilde{G}_i$ that contains $g_{\tilde{x}}$.  Let
$X_\delta = [M_i - C_\delta(M_i)] \cap \tilde{W}_i$.  Let $\bar{U}_i$ be the
closure of $\tilde{U}_i$ in $M_i$.  Then because $\partial \bar{U}_i$ is compact and
$g_{\tilde{x}}
\cap X_\delta
\neq
\emptyset$ for every $\delta > 2 \epsilon$,  we can guarantee that
$X_\delta
\subset
\tilde{U}_i$ for large enough
$\delta$.   Then $\tilde{V}_i = X_\delta$ is a
component of $M_i - C_\delta(M_i)$ which  embeds in
$N$.  Then by Lemma \ref{GFproj}, $W_i = p_i(\tilde{W}_i)$ is a component of 
$N - C_{2 \epsilon}(N)$ with boundary $p_i(\tilde{G}_i)$.  Because $\tilde{x} \in M_i
- C_{2 \epsilon}(M_i)$,
$\tilde{x} \in int \ \tilde{W}_i$.  Then
$p_i(\tilde{x})
\in int \ W_i \subset N - C_{2 \epsilon}(N)$.  But by hypothesis,
$x = p_i(\tilde{x}) \in C_{2 \epsilon}(N)$ so this is a contradiction.  Thus, 
$\tilde{U}_i \subset C_{2 \epsilon}(M_i)$ for each
$i$.

For each $i$, choose $\tilde{\mathcal T}_i$ to be a level surface in
$\tilde{U}_i$.  Recall that $\widetilde{f \circ g_i}(S_i,1)$  does not intersect
$\tilde{U}_i$ so that $\widetilde{f \circ g_i}(S_i,1)$   lies to one side of
$\tilde{\mathcal T}_i$. By construction, $\tilde{\mathcal T}_i$ is embedded, and the inclusion
map of $\tilde{\mathcal T}_i$ into $N$ is homotopic to the map $\widetilde{f \circ
g_i}: (S_i,1) \rightarrow N$. 

Let $\alpha_i^*: S_i \times [0,1] \rightarrow M_i$ be a homotopy between the inclusion
map $i_{\tilde{\mathcal T}_i}:\tilde{\mathcal T}_i
\rightarrow N$ and the map
$\widetilde{f \circ g_i}: (S_i,1) \rightarrow N$, where
$\alpha_i^*|_{(S_i, 0)} = i_{\tilde{\mathcal T}_i}$ and $\alpha_i^*|_{(S_i, 1)} =
\widetilde{f \circ g_i}|_{(S_i,1)}$. Let $\alpha_i: S_i \times [0,1] \rightarrow M_i$
be the {\it ruled homotopy\/} constructed from $\alpha_i^*$ such that for $x \in S_i$,
$\alpha_i(x,[0,1])$ is the geodesic arc with the same endpoints and in the same
homotopy class as $\alpha_i^*(x,[0,1])$. Because
$C_{2
\epsilon}(M_i)$ is convex, we know that
$\alpha(S_i,[0,1])
\subset C_{2 \epsilon}(M_i)$.  Because $\Gamma_i$ is a surface group, by Theorem
\ref{surface bound}, for $x \in \alpha_i(S_i,[0,1])$, $inj_{M_i}(x) \leq
L_{S_i} + 2 \epsilon$.

For each $i$, let ${\mathcal T}_i = p_i(\tilde{\mathcal T}_i)$ be a level surface in $U_i$. 
Consider the homotopy
$\beta_i = p_i
\circ
\alpha_i: S_i
\times [0,1]
\rightarrow N$.  Here $\beta_i|_{(S_i,0)}: S_i \rightarrow N$ is the inclusion map
$i_{{\mathcal T}_i}:{\mathcal T}_i
\rightarrow N$ which is the projection of the inclusion map $i_{\tilde{\mathcal T}_i}:
\tilde{\mathcal T}_i
\rightarrow M_i$; 
$\beta_i|_{(S_i,1)} = {f \circ g_i}|_{(S_i,1)}$; and
$\beta_i(S_i,[0,1]) \subset C_{2
\epsilon}(N)$.  By Lemma
\ref{coverinj}, for $x \in \beta_i(S_i,[0,1])$, $inj_N(x) \leq L_{S_i} + 2
\epsilon$.

Now for each $i$, choose $\tilde{T}_{i}$ to be a level surface in $\tilde{U}_i$ such
that
$(p_i|_{\tilde{U}_i})^{-1}[\bigcup_i \beta_i(S_i,[0,1])]$ lies to one side.
Let $T_i = p_i(\tilde{T}_{i})$ be a level surface in $U_i$. By construction
$(p_i|_{\tilde{U}_i})^{-1}(T_i) = \tilde{T}_i \subset \overline{C_{2
\epsilon}(M_i)}$.  

Because each of the $T_i$ is incompressible, disjoint, and homotopic to the
boundary components of
$R^*$, the \{$T_i$\} and \{$(Q_i,0)$\} span a product structure in $C_{2
\epsilon}(N)$.  (see Thm 10.5, Hempel \cite{hempel}) Hence, the \{$T_i$\} are the
boundary components of a new compact core $R$ of $N$.  In particular, since the
\{$T_i$\} lie in $C_{2 \epsilon}(N)$, we can conclude that $R \subset C_{2
\epsilon}(N)$ 

For each $i$, let $\zeta_i: S_i \times [-1,0] \rightarrow M_i$ be homotopy along the
product structure of $\tilde{U}_i$ between the inclusion maps $i_{\tilde{T}_{i}}:
\tilde{T}_{i}
\rightarrow N$ and
$i_{\tilde{\mathcal T}_i}: \tilde{\mathcal T}_i \rightarrow N$, where $\zeta_i|_{(S_i,-1)} =
i_{\tilde{T}_{i}}$ and $\zeta_i|_{(S_i,0)} =
i_{\tilde{\mathcal T}_i}$. Because $\zeta_i(S_i,[-1,0]) \subset C_{2 \epsilon}(M_i)$, by
Theorem
\ref{surface bound}, for
$x
\in
\zeta_i(S_i,[-1,0])$,
$inj_{M_i}(x) \leq L_{S_i} + 2 \epsilon$.  Then $\nu_i = p_i \circ \zeta_i: S_i
\times [-1,0]$ is a product homotopy in $U_i$ between the inclusion maps $i_{T_{i}}:
T_{i} \rightarrow N$ and $i_{{\mathcal T}_i}: {\mathcal T}_i \rightarrow N$, where
$\nu_i|_{(S_i,-1)} = i_{T_{i}}$ and $\nu_i|_{(S_i,0)} = i_{{\mathcal T}_i}$.  Furthermore,
$\nu_i(S_i,[-1,0]) \subset C_{2
\epsilon}(N)$. By Lemma \ref{coverinj}, for $x \in \nu_i(S_i,[-1,0])$, $inj_N(x)
\leq L_{S_i} + 2
\epsilon$.

Finally we will show that $$R \subset [\bigcup \beta_i(S_i \times [0,1])] \cup
[\bigcup \nu_i(S_i \times [-1,0])]$$  By Lemma \ref{what is D}, we know that $M \cong
\coprod (S_i,[-1,1])/\sim$, where
$(x_i,1) \sim (y_j,1)$ if and only if $(g_i(x_i),1) = (g_j(y_j),1)$. Let $\psi: M
\rightarrow R$ be a map such that $\psi|_{(S_i,[-1,0])} = \nu_i$ and 
$\psi|_{(S_i,[0,1])} = \beta_i$.  Then by construction, if $(g_i(x_i),1) =
(g_j(y_j),1)$, then $\psi(x_i,1) = \psi(y_j,1)$.  Hence we have a
well-defined map
$\psi: M
\rightarrow R$.
Using standard degree arguments, $\psi(M)$ is contained in the image of any proper
homotopy between $M$ and $R$.  (see Thm 2.14, Lloyd, \cite{lloyd})  In particular, $R
\subset \psi(M) \subset [\bigcup \beta_i(S_i \times [0,1])] \cup
[\bigcup \nu_i(S_i \times [-1,0])]$.

So for $x \in R$, $inj_N(x) \leq max\{L_{S_i} + 2
\epsilon\}$.  This completes the proof of Lemma \ref{R inj bdd}.
\end{proof}

Now we have bounded the injectivity radius for points in $R$.  Now consider points
in $C(N) - R$.  If $x \in C(N) - R$, then there exists an $i$, such that $x \in
U_i$.  Recall that in the previous lemma we showed that $T_i$ was a closed
incompressible separating surface and that
$(p_i|_{\tilde{U}_i})^{-1}(T_i) \subset \overline{C_{2 \epsilon}(M_i)}$. By Lemma
\ref{not in C(N)} we know that if $x \in C(N) \cap U_i$, then  $x \in
p_i(C(M_i))$.  So $inj_N(x) \leq L_{S_i}$ for some $i$.  

Then we can conclude that for $x \in C(N)$, $inj_N(x) \leq max\{L_{S_i} + 2
\epsilon\}$.  We can do this for all $\epsilon > 0$ so that for $x \in C(N)$,
$inj_N(x) \leq max\{L_{S_i}\}$.  This completes the proof of Theorem
\ref{bookresult}. 

\end{proof}

\section{The Acylindrical Case}

In this section, we will prove our main theorem in the case that $M$ is an
acylindrical, hyperbolizable 3-manifold.  Our theorem is:

\begin{theorem}
\label{acylresult}
Let $M$ be an acylindrical, hyperbolizable 3-manifold.  Then
there exists a constant $K$ such that  if $N$ is a hyperbolic 3-manifold
homeomorphic to the interior of $M$ and if $x \in C(N)$, then $inj_N(x) \leq K$.
\end{theorem}

\begin{proof}  We begin with a sketch of the
proof of the acylindrical case.  The proof will be by contradiction.  Suppose
there exists a sequence of points in the convex cores of manifolds such
that the injectivity radius based at these points goes to infinity.  We will find a
compact core $R$ in the algebraic limit which embeds in the geometric limit as
$\pi(R)$.  The compact core $\pi(R)$ will pull back
to a compact core $R_i$ in each manifold $N_i$ in the sequence such that points in
$R_i$ will have uniformly bounded injectivity radius.  The complement of
$R_i$ in $N_i$ will either be covered by the
convex cores of manifolds whose convex cores have bounded injectivity radius or
will have injectivity radius uniformly bounded by the injectivity radius of a
fixed compact subset of the geometric limit of the sequence of manifolds.  Thus,
we will have found a uniform bound on the injectivity radius for points in the
convex core of each manifold in the sequence, which is a contradiction.

Suppose for contradiction that there does not exist an upper bound on the
injectivity radius for points in the convex cores of hyperbolic 3-manifolds
homeomorphic to
$int \ M$.  Then there exists a sequence of representations \{$\rho_i: \pi_1(M)
\rightarrow Isom^+(\HHH)$\} together with its corresponding
sequence of manifolds \{$N_i = {\mathbb H}^3 / \rho_i(\pi_1(M))$\}, and a
sequence of points \{$x_i \in C(N_i)$\}, such that $\{inj_{N_i}(x_i)\}$
diverges to infinity. 

By Thurston's Compactness Theorem \ref{AH(M)cmpt}, a subsequence of 
\{$\rho_i$\} converges algebraically, up to conjugation, to a representation
$\rho$.  Let the algebraic limit manifold be $N = {\mathbb
H}^3 / \rho(\pi_1(M))$.  Using a result of Jorgensen-Marden (Prop 4.2,
\cite{jorg/marden}), we can take a further subsequence, again called
\{$\rho_i$\}, such that
\{$\rho_i(\pi_1(M))$\} converges to $\hat{\Gamma}$ geometrically. Let the
geometric limit manifold be
$\hat{N} = {\mathbb H}^3 /
\hat{\Gamma}$.  By definition of geometric convergence, there exists a sequence
of $(K_i,r_i)$-approximate isometries $f_i: B_{r_i}(0) \subset N_i
\rightarrow \hat{N}$ such that  $K_i \rightarrow 1$ and $r_i \rightarrow
\infty$. Furthermore, because
$\rho(\pi_1(M))$ is a subgroup of $\hat{\Gamma}$, there is a natural covering map from
the algebraic limit to the geometric limit
$\pi: N \rightarrow \hat{N}$.

Let $N^* = (\HHH \cup \Omega_{\rho(\pi_1(M)})/ \rho(\pi_1(M))$ be the conformal
extension of $N$.  Anderson and Canary
\cite{andcanary} have an alternate definition of an accidental parabolic which,
to avoid confusion, we will call an {\it unexpected parabolic\/}.
We say that $\rho(\pi_1(M))$ has connected limit set and no
unexpected parabolics if and only if every closed curve $\gamma$ in
$\partial N^*$ which is homotopic to a curve of arbitrarily small length in
$N^*$ is homotopic to a curve of arbitrarily small length in $\partial
N^*$.  

The following lemma shows that the fundamental group of an acylindrical,
hyperbolizable manifold has connected limit set and no unexpected parabolics. 
These properties will be useful in showing the existence of a compact core in the
algebraic limit which will embed in the geometric limit.

\begin{lemma}
\label{acylnopara}  Let $M$ be an acylindrical, hyperbolizable
3-manifold.  Let $\rho \in {\mathcal D}(\pi_1(M))$, $N = \HHH
/ \rho(\pi_1(M))$, and $N^* = (\HHH \cup
\Omega_{\rho(\pi_1(M))})/ \rho(\pi_1(M))$. Then every closed curve $\gamma$ in
$\partial N^*$ which is homotopic to a curve of arbitrarily small length in
$N^*$ is homotopic to a curve of arbitrarily small length in $\partial N^*$. 
Therefore, $\rho(\pi_1(M))$ has connected limit set and no unexpected
parabolics.
\end{lemma}

\begin{proof}  The proof will be by
contradiction.  For $0 < \epsilon < \epsilon_3$, McCullough
\cite{mccullough} and Kulkarni-Shalen \cite{kul/shalen} guarantee the existence of a
relative compact core
$R^*$ of
$\overline{C^\circ_\epsilon(N)}$ with associated parabolic locus $P = \partial R^* \cap
\partial N_{thin(\epsilon)}$ such that
$\partial C(N) \cap C^\circ_\epsilon(N) \subset \partial R^*$.  Note that
$\pi_1(N) = \pi_1(\overline{C_{\epsilon}(N)})$, so that $R^*$ is also a compact
core of $N$.  Because $N$ is homeomorphic to the interior of $M$ and $R^*$ is a
compact core of $N$, we can conclude that
$R^*$ is homeomorphic to $M$.  (Thm 1, McCullough-Miller-Swarup
\cite{mmswarup}) Therefore, because $M$ has incompressible boundary and is
acylindrical, 
$R^*$ has incompressible boundary and is acylindrical. 

%already showed that $M$ and $R^*$ are homeo
%Because $R^*$ is an orientable, irreducible, compact manifold with boundary that
%is homotopy equivalent to $M$, the following theorem due to Johannson shows that
%$R^*$ is homeomorphic to $M$.  Hence, $R^*$ is also acylindrical.
%\begin{thm}
%\label{acylhomeo} (Corollary to Thm 24.2, Johannson \cite{johannson}) Let $M$ and
%$M^\prime$ be orientable, irreducible, compact 3-manifolds with boundary.  Suppose that
%$h: M
%\rightarrow M^\prime$ is a map such that $h_*: \pi_1(M) \rightarrow
%\pi_1(M^\prime)$ is an isomorphism.  If $M$ is acylindrical, then there is a
%homotopy $h_t: M
%\rightarrow M^\prime$ such that $h_0 = h$ and $h_1: M
%\rightarrow M^\prime$ is a homeomorphism.
%\end{thm}

Suppose for contradiction that there exists a homotopically non-trivial, closed
curve
$\gamma$ in
$\partial N^*$ which is homotopic to a curve of arbitrarily
small length in $N^*$, but that $\gamma$ is not homotopic to a curve of
arbitrarily small length in $\partial N^*$. 

Recall that Sullivan \cite{sullivan2} has shown that there exists a homeomorphism
$g: N^*
\rightarrow C(N)$ which is homotopic to the canonical nearest point retraction
$\hat{r}: N^* \rightarrow C(N)$ (see Sec 1.3, Epstein-Marden
\cite{epstein/marden}), and which is
$K$-bilipschitz on $\partial N^*$ where $K$ is independent of
$\Gamma$. 

Suppose $\gamma$ is homotopically trivial in $N^*$, then $g(\gamma)$ is a
homotopically non-trivial, closed curve in $\partial C(N)$ which is homotopically
trivial in $C(N)$.  Here, $\partial C(N) \cap C^\circ_\epsilon(N)$ is a compact
core of
$\partial C(N)$ so that there exists a homotopically non-trivial curve in 
$\partial C(N) \cap C^\circ_\epsilon(N)$ which is homotopically trivial in $C(N)$.  But
$\partial C(N) \cap C^\circ_\epsilon(N)$ is an incompressible subsurface of $\partial
R^*$, and  $R^*$ has incompressible boundary, so we can conclude that $\partial C(N)
\cap C^\circ_\epsilon(N)$ is incompressible in
$C(N)$. But this is a contradiction.

Therefore $\gamma$ is homotopically non-trivial in $N^*$.  Because $g$ is
$K$-bilipschitz on $\partial N^*$, we know that $g(\gamma)$ is a closed curve in
$\partial C(N)$ which is homotopic to a curve of arbitrarily small length in
$C(N)$, but that $g(\gamma)$ is not homotopic to a curve of arbitrarily small
length in $\partial C(N)$.  Because $\partial C(N) \cap C^\circ_\epsilon(N)$ is a compact
core of
$\partial C(N)$, we know that $g(\gamma)$ is homotopic to a curve on $\partial C(N)
\cap C^\circ_\epsilon(N) \subset \partial R^*$.  Therefore, without loss of
generality, we can consider $g(\gamma)$ to be a closed curve on $\partial R^*$.
Because
$g(\gamma)$ is homotopic to a curve of arbitrarily small length in $C(N)$, there
exists a closed curve $\alpha \subset P$ such that
$g(\gamma)$ is homotopic to $\alpha$ in $C(N)$.  Then since $P \subset \partial
R^*$ and $R^*$ is a compact core, we can conclude that
$g(\gamma)$ is homotopic to $\alpha$ in $R^*$.

Suppose $g(\gamma)$ is homotopic to $\alpha$ in $\partial
R^*$.  Then $g(\gamma)$ is homotopic into the parabolic locus
of $R^*$ in $\partial R^*$.  Then $g(\gamma)$ is a peripheral curve in
$\partial C(N)$, and hence is homotopic to a curve of arbitrarily small length
in $\partial C(N)$. Then
$\gamma$ must have been homotopic to a curve of arbitrarily small length in
$\partial N^*$.  But this contradicts our initial assumptions about
$\gamma$.   So $g(\gamma)$ is not homotopic to $\alpha$ in $\partial R^*$.

Using the Homotopy Annulus Theorem (Thm VIII.10, Jaco
\cite{jaco}), we can construct a $\pi_1$-injective, proper embedding $f: S^1 \times
[0,1] \rightarrow R^*$ such that $f(S^1,0) = g(\gamma)$ and $f(S^1,1) = \beta$ where 
$f(S^1,[0,1])$ cannot be properly homotoped into $\partial R^*$.  Therefore, $R^*$
is not acylindrical. But this is also a contradiction. Thus, every closed curve
$\gamma$ in $\partial N^*$ which is homotopic to a curve of arbitrarily small
length in
$N^*$ is homotopic to a curve of arbitrarily small length in $\partial N^*$. 
Therefore, $\rho(\pi_1(M))$ has connected limit set and no unexpected
parabolics. This completes the proof of Lemma
\ref{acylnopara}.
\end{proof}

Next, we will show that there exists a compact core $R$ in the algebraic limit
$N = \HHH /
\rho(\Gamma)$ which embeds in the geometric limit $\hat{N} = \HHH / \hat{\Gamma}$
as $\pi(R)$.

\begin{lemma}
\label{ACcore} Let $M$ be an acylindrical, hyperbolizable
3-manifold.  Let $\pi: N \rightarrow \hat{N}$ be the covering map between the
algebraic limit $N = \HHH / \rho(\pi_1(M))$ and the geometric limit
$\hat{N} = \HHH / \hat{\Gamma}$. Then there exists a compact core $R \subset N$
such that $\pi(R)$ embeds in $\hat{N}$.
\end{lemma}

\begin{proof} There are two cases---either the
limit set of $\rho(\pi_1(M))$ is the entire sphere or not.

Suppose the limit set of $\rho(\pi_1(M))$ is the entire sphere.  Because
$\rho(\pi_1(M))$ satisfies Bonahon's Condition (B), 
$N = \HHH / \rho(\pi_1(M))$ is topologically tame.  
Then the following theorem of Canary states that the algebraic limit and the
geometric limit agree, and, in this case, any compact core of
the algebraic limit is also a compact core of geometric limit by default.

\begin{theorem}
\label{limit set all} (Thm 9.2, Canary \cite{canary1})  Let $\{\rho_i:
\pi_1(M) \rightarrow Isom^+(\HHH)\}$ be a sequence of discrete faithful
representations converging algebraically to $\rho: \pi_1(M) \rightarrow
Isom^+(\HHH)$.  If the limit set of $\rho(\pi_1(M))$ is all of $S^2_\infty$ and $N =
\HHH / \rho(\pi_1(M))$ is topologically tame, then $\{\rho_i\}$ converges strongly to
$\rho$.
\end{theorem}

If the limit set of $\rho(\pi_1(M))$ is not the entire sphere, then the domain
of discontinuity of $\rho(\pi_1(M))$ has nonempty domain of discontinuity. 
Recall that by Lemma \ref{acylnopara}, $\rho(\pi_1(M))$ has connected limit set and
no unexpected parabolics.  In this case, the following theorem of Anderson-Canary
guarantees that given an algebraically convergent sequence such that its
associated image groups converge geometrically, we can find a compact core in the
algebraic limit which embeds in the geometric limit.

\begin{theorem}
\label{limit set some}  (Cor B, Anderson-Canary \cite{andcanary})  Let
$\pi_1(M)$ be a finitely generated, torsion-free, nonabelian group, and let
$\{\rho_i\}$ be a sequence in ${\mathcal D}(\pi_1(M))$ converging algebraically to
$\rho$.  Suppose that
$\{\rho_i(\pi_1(M))\}$ converges geometrically to $\hat{\Gamma}$.  Let $N = \HHH
/ \rho(\pi_1(M))$, $\hat{N} = \HHH / \hat{\Gamma}$, and let $\pi: N \rightarrow
\hat{N}$ be the covering map.  If $\rho(\pi_1(M))$ has nonempty domain of
discontinuity, connected limit set, and contains no unexpected parabolics,
then there exists a compact core $R$ of $N$ such that $\pi: N \rightarrow \hat{N}$ is an
embedding restricted to $R$.
\end{theorem}

Thus, in either case, we can find a compact core in the algebraic limit which
embeds in the geometric limit.  This completes the proof of Lemma
\ref{ACcore}.
\end{proof}

By Lemma \ref{acylnopara}, we know that $\rho(\pi_1(M))$ has connected limit set
and contains no unexpected parabolics.  The next lemma due to Canary-Minsky and
Anderson-Canary shows that for large enough $i$, the pull back of $\pi(R)$ to
$N_i$ is a compact core of $N_i$.

\begin{lemma} 
\label{core in N_i} (Lem 7.2, Anderson-Canary \cite{andcanary})  Let $\pi_1(M)$
be a finitely generated, torsion-free, nonabelian group, and let $\{\rho_i\}$
be a sequence in ${\mathcal D}(\pi_1(M))$ converging algebraically to $\rho$. 
Suppose that $\{\rho_i(\pi_1(M))\}$ converges geometrically to $\hat{\Gamma}$. 
Let $N = \HHH / \rho(\pi_1(M))$, $\hat{N} = \HHH / \hat{\Gamma}$, and let $\pi: N
\rightarrow
\hat{N}$ be the covering map.  Suppose $\rho(\pi_1(M))$ has connected limit set
and contains no unexpected parabolics.  Let
$R$ be a compact core of $N$ such that $\pi$ is an embedding restricted to $R$.
Then for large enough $i$, $R_i = f_i^{-1}(\pi(R))$ is a compact core for $N_i$.
\end{lemma}

Now let us show that points in $R_i$ have uniformly
bounded injectivity radius, where the bound depends on the compact set $\pi(R)$ in
$\hat{N}$.

\begin{lemma}
\label{injbddinR_i}  For large enough $i$, and for $x \in R_i =
f_i^{-1}(\pi(R)) \subset N_i$, we have $inj_{N_i}(x) \leq 2 \kappa_R$.
\end{lemma}

\begin{proof}  Note that because $\pi(R)$ is a
compact set in $\hat{N}$, there exists a constant $\kappa_R$ such that for $x \in
\pi(R)$, we have $inj_{\hat{N}}(x) \leq \kappa_R$. Then, because $R$ is
compact and
$f_i: B_{r_i}(0) \rightarrow \hat{N}$ is a $(K_i,r_i)$-approximate isometry
such that $K_i \rightarrow 1$ and $r_i \rightarrow \infty$, it is
possible to choose $I > 0$ such that for $i > I$, we guarantee that $K_i < 2$ and
that the closure of the
$4\kappa_R$-neighborhood of $R_i$ lies in $B_{r_i}(0)$.

Let  $\gamma_x$ be a homotopically non-trivial loop
in $\hat{N}$ that is based at $x \in \pi(R)$ and is of length
$\leq 2
\kappa_R$.  If
$f_i^{-1}(\gamma_x)$ is a homotopically trivial loop, then
$f_i^{-1}(\gamma_x)$ bounds a disk $D_x$ in
$N_i$.  Consider the immersion $g: D^2 \rightarrow D_x$.  For $y,z
\in
\partial D^2$, let $yz$ be the line segment in $D^2$ joining $y$ and $z$.  Let
$g(yz)^*$ be the geodesic arc in $N_i$ that is properly homotopic to the segment
$g(yz)$.  Fix $z \in D^2$.  Then the new disk $D_x^\prime =
\bigcup_y g(yz)^*$ has diameter $\leq 2 K_i \kappa_R$ in $N_i$.  For $i > I$, we
know that $K_i < 2$ so that
$D_x^\prime \subset B_{r_i}(0)$.  Then $f_i(D_x^\prime)$ is
a disk in $\hat{N}$, so $\gamma_x$ is homotopically trivial in
$\hat{N}$.  But this is a contradiction. Therefore, $f_i^{-1}(\gamma_x)$ is a
homotopically non-trivial loop based at $f_i^{-1}(x)$ of length $\leq 2 K_i
\kappa_R$.  Then for $x \in R_i$, we know that $inj_{N_i}(x) \leq K_i \kappa_R$. 
In particular, for $i > I$, we know that $K_i < 2$, so for $x \in R_i$,
we have $inj_{N_i}(x)
\leq 2
\kappa_R$.  
\end{proof}

The following lemma shows that points in $C(N_i) - R_i$ also have bounded
injectivity radius, because they are either covered by manifolds whose convex cores
have bounded injectivity radius or have injectivity radius bounded by the
injectivity radius of a fixed compact subset of the geometric limit.

\begin{lemma}
\label{not in R_i} Let $0 < \epsilon < \epsilon_3$, and let \{$S_j$\} be the
boundary components of $\pi(R) \subset \hat{N}$.  Then for large enough $i$ and
$x
\in C(N_i) - R_i$, we have $inj_{N_i}(x)
\leq max \{L_{S_j}, 2 \kappa_{S_j}, \epsilon\}$.
\end{lemma}

\begin{proof} Temporarily fix $i$. Because $N_i$ is
homeomorphic to the interior of $M$ where $M$ has incompressible boundary, the
compact core
$R_i$ of $N_i$ is homeomorphic to $M$.  (Thm 1, McCullough-Miller-Swarup
\cite{mmswarup}) Since $M$ has incompressible boundary, so does
$R_i$.  Let
\{$T_{ij}$\} be the components of $\partial R_i$.  Then each $T_{ij} =
f_i^{-1}(S_j)$ is homeomorphic to $S_j$ and is an incompressible separating
surface of $N_i$.

Let $U_j$ be the component of $N_i - int \ R_i$
with boundary component $T_{ij}$.  Because $R_i$ is a compact core of a
topologically tame 3-manifold with incompressible boundary, by Lemma
\ref{prodstructure}, $U_j$ possesses a product structure $T_{ij} \times
[0,\infty)$ for each $j$.  Because $i_*(\pi_1(T_{ij})) = \pi_1(M_j)$, by the
Lifting Theorem, the inclusion map $i: T_{ij}
\times [0,\infty)
\rightarrow U_j$ lifts to a map $\tilde{i}: T_{ij} \times [0,\infty) \rightarrow
M_j$. Let $\tilde{U}_j = \tilde{i}(T_{ij}, [0,\infty))$.  Then the
projection map $p_j|_{\tilde{U}_j}: \tilde{U}_j \rightarrow U_j$ is a
homeomorphism.  Let $\tilde{T}_{ij} = (p_j|_{\tilde{U}_j})^{-1}(T_{ij})$.  

If $x \in C(N_i) - R_i$, then there exists $j$ such that $x \in U_j$. 
By Lemma \ref{coverinj}, it suffices to bound the injectivity radius based at
$\tilde{x} = (p_j|_{\tilde{U}_j})^{-1}(x)$ in
$M_j$.  There are two possibilities: \ either $\tilde{x} \in C(M_j)$ or not.
If $\tilde{x} \in C(M_j)$, then by Theorem \ref{surface bound},
$inj_{N_i}(x) \leq L_{S_j}$.  

Suppose $\tilde{x} \in M_j - C(M_j)$.  There are two possibilities here also: \
either
$\tilde{x}
\in (M_j)_{thin(\epsilon)}$ or $\tilde{x} \in (M_j)^\circ_\epsilon - C(M_j)$. 
If $\tilde{x} \in (M_j)_{thin(\epsilon)}$, then
$inj_{M_j}(\tilde{x}) \leq \epsilon$.  

Suppose $\tilde{x}
\in (M_j)^\circ_\epsilon - C(M_j)$.  Then there exists a geodesic ray
$\bar{g}_{\tilde{x}}$ that is perpendicular to a component $\tilde{A}_j$ of $\partial
C(M_j)$ and that passes through $\tilde{x}$.  Let $g_{\tilde{x}}$ be the portion of
$\bar{g}_{\tilde{x}}$ beginning at $\tilde{x}$. Note that by Lemma \ref{expinj}, the
injectivity radius strictly increases out a geometrically finite end.  Because
$\tilde{x} \in (M_j)^\circ_\epsilon$, the ray
$g_{\tilde{x}}$ is entirely contained in $(M_j)^\circ_\epsilon$.  Then either
$g_{\tilde{x}}$ is contained in
$\tilde{U}_j$, or $g_{\tilde{x}}$ intersects $\tilde{T}_{ij}$.

Suppose $g_{\tilde{x}}$ is contained in
$\tilde{U}_j$.  Let
$\tilde{W}_j$ be the closure of the component of
$M_j - \tilde{A}_j$ that contains $g_{\tilde{x}}$.  Let
$X_\delta = [M_j - C_\delta(M_j)] \cap \tilde{W}_j$.  Because $\tilde{T}_{ij}$ is
compact and
$g_{\tilde{x}} \cap X_\delta \neq \emptyset$ for all $\delta > 0$, we can guarantee
that
$X_\delta \subset
\tilde{U}_j$ for large enough $\delta$.  Then $\tilde{V}_j = X_\delta$ is a
component of $M_j - C_\delta(M_j)$ that embeds in
$N_i$.  By Lemma \ref{GFproj}, $W_j = p_j(\tilde{W}_j)$ is the component of $N_i -
C(N_i)$ with boundary $p_j(\tilde{A}_j)$.  Because $\tilde{x} \in M_j - C(M_j)$,
$\tilde{x} \in int \ \tilde{W}_j$.  Then
$p_j(\tilde{x})
\in int \ W_j \subset N_i - C(N_i)$.  But by hypothesis,
$x = p_j(\tilde{x}) \in C(N_i)$ so this is a contradiction.  

Thus $g_{\tilde{x}}$ must intersect $\tilde{T}_{ij}$.  By Lemma
\ref{expinj}, we know that the injectivity radius strictly increases out a
geometrically finite end. Thus, it suffices to bound the injectivity radius for
points $\tilde{y} \in
\tilde{T}_{ij}$.
  
Because $S_j$ is in the image of a compact set in $\hat{N}$, there exists a
constant $\kappa_{S_j}$ such that for $z \in S_j$, $inj_{S_j}(z) \leq
\kappa_{S_j}$.  Because $T_{ij} = f_i^{-1}(S_j)$ and $f_i$ is a
$(K_i,r_i)$-approximate isometry, by the argument in the proof of Lemma
\ref{injbddinR_i}, for
$y
\in T_{ij}$, we know 
$inj_{T_{ij}}(y)
\leq K_i \kappa_{S_j}$.  Furthermore, because $p_j|_{\tilde{U}_j}$ is an isometry
and
$\tilde{T}_{ij} =  (p|_{\tilde{U}_j})^{-1}(T_{ij})$, for
$\tilde{y} \in \tilde{T}_{ij}$, we know
$inj_{\tilde{T}_{ij}}(\tilde{y})
\leq K_i
\kappa_{S_j}$. 
Because $T_{ij}$ is incompressible in $N_i$ and $\tilde{T}_{ij}$ is a lift of $T_{ij}$,
$\tilde{T}_{ij}$ is incompressible in $M_j$.  So for
$\tilde{y} \in \tilde{T}_{ij}$, we know that
$inj_{M_j}(\tilde{y})
\leq inj_{\tilde{T}_{ij}}(\tilde{y})$.  

Then by Lemmas
\ref{coverinj} and \ref{expinj}, for $x \in C(N_i) - R_i$, $\tilde{x}
\in (M_j)^\circ_\epsilon - C(M_j)$, and $g_{\tilde{x}}$ intersecting
$\tilde{T}_{ij}$ at $\tilde{y}$, we have
$inj_{N_i}(x) \leq inj_{M_j}(\tilde{x}) \leq inj_{M_j}(\tilde{y}) \leq K_i
\kappa_{S_j}$.  Because $K_i \rightarrow 1$, for large $i$, we have $inj_{N_i}(x)
\leq 2 \kappa_{S_j}$.

Thus, for large $i$ and $x \in C(N_i) - R_i$, we can conclude that $inj_{N_i}(x)
\leq max \{L_{S_j}, 2 \kappa_{S_j}, \epsilon\}$.  This completes the proof of
Lemma \ref{not in R_i}. 
\end{proof}

Therefore, by Lemmas \ref{injbddinR_i} and \ref{not in R_i}, for $x
\in C(N_i)$, we know that $inj_{N_i}(x) \leq max\{2\kappa_R, L_{S_j}, 2
\kappa_{S_j},
\epsilon\}$.  We can do this for all $0 < \epsilon < \epsilon_3$, so that for $x
\in C(N_i)$, we have
$inj_{N_i}(x) \leq max\{2\kappa_R, L_{S_j}, 2 \kappa_{S_j}\}$.  This uniform
bound contradicts the assumption that there exists a sequence of points \{$x_i \in
C(N_i)$\} such that
\{$inj_{N_i}(x_i)$\} converges to infinity.  This completes the proof of Theorem
\ref{acylresult}. 
\end{proof}

\section{Some Consequences}

In this section, we present some consequences of the main theorem. 
First we present a slightly stronger answer to  McMullen's conjecture in the
case of a book of $I$-bundles.  We will show that if $N$ is a hyperbolic
3-manifold homotopy equivalent to a book of $I$-bundles, then there exists an
upper bound on injectivity radius for points in the convex core of
$N$, where the bound depends on the number of generators in $\pi_1(N)$. 

\begin{corollary}
\label{conjresult}
Let $N$ be a hyperbolic 3-manifold homotopy equivalent to a book of
$I$-bundles.   Then there exists a constant $L^\prime$ such that for $x
\in C(N)$,
$inj_N(x) \leq L^\prime$, where $L^\prime$ depends on the number of
generators of $\pi_1(N)$.
\end{corollary}

\begin{proof} First let us show 
that if $N$ is homotopy equivalent to a book of $I$-bundles, then $N$ is
homeomorphic to the interior of a book of $I$-bundles.  Let $R$ be a compact
core for
$N$.  Because $\pi_1(N)$ satisfies Bonahon's Condition (B), $N$ is
topologically tame, and hence $N$ is homeomorphic to the interior of its
compact core $R$.  Thus, it suffices to show that $R$ is a book of $I$-bundles. 

We will use the characteristic submanifold theory developed by Johannson
\cite{johannson} and Jaco-Shalen \cite{jaco/shalen}.  First let us
introduce some definitions.  A map
$f:(V,\partial V) \rightarrow (M, \partial M)$ of an annulus or torus, $V$,
into $M$ is {\it essential\/} if $f$  is
$\pi_1$-injective and $f(V)$ is not properly homotopic into $\partial M$.  
A map $f:E \rightarrow M$ of an $I$-bundle $E$ into $M$ is {\it
admissible\/} if $f^{-1}(\partial M)$ is the associated $\partial I$-bundle. 
A map $f: E \rightarrow M$ of an $I$-bundle $E$ into $M$ is {\it essential\/} if
$f$ is
$\pi_1$-injective, and for every component $V$ of $\partial E -
f^{-1}(\partial M)$, $f|_{V}$ is an essential map of an annulus or torus into
$M$. 

A compact submanifold $\Sigma$ of $M$ is a {\it characteristic submanifold\/}
if $\Sigma$ is a minimal collection of admissibly embedded, essential
$I$-bundles and Seifert fiber spaces such that every essential,
admissible embedding $f: E \rightarrow M$ of a Seifert fiber space or
$I$-bundle into $M$ is properly homotopic to an admissible map with image in
$\Sigma$.  In fact, every compact, orientable, irreducible 3-manifold with
incompressible boundary contains  a unique (up to isotopy) characteristic
submanifold.  

We can apply the characteristic submanifold theory to $M$ and to $R$ to obtain
the characteristic submanifolds $\Sigma_M$ and $\Sigma_R$.   
Recall that $M$ is a  book of $I$-bundles if there exists a disjoint collection
$A$ of incompressible annuli such that each component of the manifold obtained
by cutting $M$ along $A$ is either a solid torus, or an $I$-bundle $R$ over a
surface of negative Euler characteristic such that $\partial R \cap \partial M$
is the associated $\partial I$-bundle.  Thus, we know that $\Sigma_M$ is a
collection of $I$-bundles and solid tori, and $M - \Sigma_M$ is a collection of
solid tori.

Because $R$ is a hyperbolizable 3-manifold and $\pi_1(R)$ contains no $\Zz
\oplus \Zz$ subgroup, then (Sec 11, Morgan \cite{morgan}) $\Sigma_R$ is a
collection of $I$-bundles and solid tori.  Let $H_0: M \rightarrow R$ be a
homotopy equivalence.  By Johannson (Thm 24.2, \cite{johannson}), the map $H_0: M
\rightarrow R$ is homotopic to a map $H_1: M
\rightarrow R$ such that $H_1^{-1}(\Sigma_R) = \Sigma_M$, $H_1|_{M - \Sigma_M}:
M -
\Sigma_M
\rightarrow R - \Sigma_R$ is a homeomorphism, and $H_1|_{\Sigma_M}: \Sigma_M
\rightarrow \Sigma_R$ is a homotopy equivalence.  Because $H_1:
M - \Sigma_M \rightarrow R - \Sigma_R$ is a homeomorphism,
each component of
$R - \Sigma_R$ is a solid torus.  Then as a consequence of a result of
Culler-Shalen (Prop 4.3, \cite{culler/shalen}), $R$ is also a book of
$I$-bundles.  Thus, $N$ is homeomorphic to the interior of a book of $I$-bundles
$R$. 

Now we will show a relationship between the Euler characteristic of the
boundary components of $R$ and the number of generators of $\pi_1(R)$.  Let
$DR$ be the double of $R$.  Then
$DR$ is a closed 3-manifold which has Euler characteristic $\chi(DR) = 0$.  So
$\chi(DR) = 2 \chi(R) - \chi(\partial R)$ or $\chi(\partial R) = 2 \chi(R)$. 
Recall $\chi(R) = \beta_0 - \beta_1 + \beta_2 - \beta_3$ where $\beta_i$ is the
 rank of $H_i(R)$.  Because $R$ is a
connected 3-manifold with boundary, $\chi(R) = 1 - \beta_1 + \beta_2 \geq
1 - \beta_1$.  Because $H_1(R)$ is the abelianization of
$\pi_1(R)$, $\beta_1 \leq n$ where $n$ is the number of generators in
$\pi_1(R)$.  So
$\chi(R) \geq 1 - n$ or $\chi(\partial R) \geq 2 - 2n$.  

By applying Theorem \ref{bookresult} to $R$ and $N$, we obtain an upper
bound $L$ such that for $x \in C(N)$, $inj_N(x) \leq L$.  Recall from the proof
of Theorem \ref{bookresult}, $L = max \{L_{S_j}\}$, where the maximum is taken
over all boundary components $\{S_j\}$ of $R$ and
$L_{S_j}$ is the bound obtained in Theorem \ref{surface bound}.
Let $L^\prime = max \{L_{S_j}\}$ where the maximum is taken over all
surfaces $\{S_j\}$ such that $|\chi(S_j)| \leq 2n - 2$.  Note that this
is a finite set of surfaces.  Then $L \leq L^\prime$, and hence for $x
\in C(N)$,
$inj_N(x) \leq L^\prime$, where the bound $L^\prime$ depends only on the
number of generators of $\pi_1(R) =
\pi_1(N)$.  This completes the proof of Corollary
\ref{conjresult}.
\end{proof}

\begin{remark} 
Note that McMullen's conjecture concerns the radius of balls
embedded in $C(N)$ rather than injectivity radius which involves the radius of 
balls embedded in $N$.  We can see the necessity for this by considering the case
when $M$ is a handlebody of genus 2.  Let $\Gamma_i$ be a free group on two
generators constructed as follows: \ in the ball model of hyperbolic 3-space,
let
$\Gamma_i$ be generated by hyperbolic isometries which identify two pairs of
disjoint hemispheres which are perpendicular to
$S^2_\infty$ and whose fixed points are antipodally situated on $S^2_\infty$. 
Then $N_i = \HHH /
\Gamma_i$ is a hyperbolic 3-manifold homotopy equivalent to $M$. In this case,
a fundamental domain in
$\HHH$ for the action of $\Gamma_i$ is the portion of $\HHH$ lying ``outside''
the hemispheres.  Because the fixed points are antipodally situated on
$S^2_\infty$, the origin is in
$CH(\Lambda_{\Gamma_i})$ and hence its projection $p_i(0)$ will lie in
$C(N_i)$.  The injectivity radius based at $p_i(0)$ in
$N_i$ is greater than or equal to the radius of the largest ball in $\HHH$ based
at the origin that can be embedded in the fundamental domain.  As $i \rightarrow
\infty$, let the Euclidean radius of the hemispheres shrink to 0.  Then as $i
\rightarrow \infty$, 
$inj_{N_i}(p_i(0)) \rightarrow \infty$ so that a uniform upper bound does
not exist over all hyperbolic 3-manifolds homotopy equivalent to $M$.

When $i$ is large, however, $CH(\Lambda_{\Gamma_i})$ is very ``thin and long''
so that large balls based at the origin cannot be embedded in
$CH(\Lambda_{\Gamma_i})$, and hence large balls based at the origin cannot be
embedded in $C(N)$.  Because of this example, McMullen only considered the
radius of balls embedded in $C(N)$, rather than the injectivity radius which
involves the radius of balls embedded in $N$.  
 
Despite this example, in some cases it is still possible to find a uniform
upper bound on injectivity radius.  Note that in the example which failed
to have uniformly bounded injectivity radius, $M$ had compressible boundary. In
this paper, we only considered the case that
$M$ had incompressible boundary. 
\end{remark}

The Main Theorem, along with a result of McMullen, shows that the
limit set varies continuously over the space of hyperbolic 3-manifolds of the
same topological type.  

\begin{corollary}
\label{varylimset1} Let $M$ be a book of I-bundles or an
acylindrical, hyperbolizable 3-manifold.  Let \{$N_i = {\mathbb H^3} /
\Gamma_i$\} be a sequence of hyperbolic 3-manifolds with base frame $\omega_i$ in
$C(N_i)$ such that each $N_i$ is homeomorphic to the interior of
$M$ and the injectivity radius at $\omega_i$ is bounded away from $0$.  If
\{$N_i$\} converges geometrically to
$N = {\mathbb H^3} /
\Gamma$, then
\{$\Lambda_{\Gamma_i}$\} converges to $\Lambda_\Gamma$ in the Hausdorff
topology.
\end{corollary}

\begin{proof} The proof is a direct corollary of the Main Theorem and Prop 2.4,
McMullen
\cite{mcmullen} which relates upper and lower bounds on injectivity radius
in the convex core to  convergence of limit sets.
\end{proof}

%\begin{theorem}
%\label{varylimset} (Prop 2.4, McMullen \cite{mcmullen}) 
%For $0<r<R$, let \{$N_i =  \HHH / \Gamma_i$\} be a sequence of
%hyperbolic 3-manifolds with base frame $e_i$ such that:
%\begin{enumerate}
%\item the baseframe is based in $C(N_i)$,
%\item the injectivity radius at $e_i$ is greater than $r$, and
%\item for $x \in C(N_i)$, we know $inj_{N_i}(x) \leq R$.
%\end{enumerate} Suppose \{$N_i$\} converges geometrically to a
%limit manifold $N =  \HHH / \Gamma$.   Then \{$\Lambda_{\Gamma_i}$\}
%converges to $\Lambda_\Gamma$ in the Hausdorff topology.

% For $0<r<R$,
%let
%${\mathcal H}^n_{r,R}$ denote the space of hyperbolic $n$-manifolds with base
%frame
%$(M,e)$ such that:
%\begin{enumerate}
%\item the baseframe $e$ is in the convex core of $M$,
%\item at $e$, the injectivity radius in $M$ is greater than $r$, and
%\item for points in the convex core, the injectivity radius in $M$ is bounded
%above by $R$.
%\end{enumerate}  Let \{$N_i = {\mathbb H}^n / \Gamma_i$\} be a sequence of
%$n$-manifolds in ${\mathcal H}^n_{r,R}$ that converges geometrically to a limit
%manifold $N = {\mathbb H}^n / \Gamma$.   Then \{$\Lambda_{\Gamma_i}$\}
%converges to $\Lambda_\Gamma$ in the Hausdorff topology. 
%\end{theorem}
%
%\begin{proof} Let us sketch the proof of this theorem.  Let $T(\Gamma,R) = \{x
%\in \HHH: inj_N(x) \leq R\}$.  Then $\overline{T(\Gamma,R)} \cap S^2_\infty
%\subset \Lambda_\Gamma$ because each point in $T(\Gamma,R)$ is close to
%$C(N)$ or $N_{thin(\epsilon)}$ and a sequence of points in $T(\Gamma,R)$ which
%limits to $S^2_\infty$ must limit to a point in $\Lambda_\Gamma$.
%
%\end{proof}

Now we will present another corollary of the main theorem which does not involve
base frame considerations:

\begin{corollary}
\label{limitset2}
Let $M$ be a book of I-bundles or an
acylindrical, hyperbolizable 3-manifold.  Let \{$N_i = {\mathbb H^3} /
\Gamma_i$\} be a sequence of hyperbolic 3-manifolds homeomorphic to the
interior of
$M$.  If \{$N_i$\} converges geometrically to $N = {\mathbb H^3} / \Gamma$ and
$\Gamma$ is nonabelian, then \{$\Lambda_{\Gamma_i}$\} converges to
$\Lambda_\Gamma$ in the Hausdorff topology.
\end{corollary}

\begin{proof}  
Let us first state and outline the proof of the result of McMullen cited in the
proof of the previous theorem. 

\begin{theorem}
\label{varylimset} (Prop 2.4, McMullen \cite{mcmullen}) 
For $0<r<R$, let \{$N_i =  \HHH / \Gamma_i$\} be a sequence of
hyperbolic 3-manifolds with base frame $\omega_i$ such that:
\begin{enumerate}
\item the baseframe $\omega_i$ lies in $C(N_i)$,
\item the injectivity radius at $\omega_i$
is greater than $r$, and
\item for $x \in C(N_i)$, we require that $inj_{N_i}(x)$ be bounded above by $R$.
\end{enumerate} Suppose \{$N_i$\} converges geometrically to a
limit manifold $N =  \HHH / \Gamma$.   Then \{$\Lambda_{\Gamma_i}$\}
converges to $\Lambda_\Gamma$ in the Hausdorff topology.
\end{theorem}

\begin{proof} Let
$T(\Gamma,R) = \{x \in \HHH: inj_N(x) \leq R\}$.  Because each point in
$T(\Gamma,R)$ is a bounded distance (depending only on R) away from $C(N)$ or
$N_{thin(\epsilon)}$, we know that a limit point of $T(\Gamma_i,R)$ in
$S^2_\infty$ must be a point in
$\Lambda_\Gamma$, that is, $\overline{T(\Gamma,R)}
\cap S^2_\infty \subset \Lambda_\Gamma$. 

By geometric convergence, $inj_{N_i}(x)$ converges uniformly
to $inj_N(x)$ on compact subsets of $\HHH$.  Therefore $\lim \sup T(\Gamma_i,R)
\subset T(\Gamma,R)$.  

Because injectivity radius for points in the convex core is
bounded above by $R$, we know that $CH(\Lambda_{\Gamma_i}) \subset
T(\Gamma_i,R)$.
Without loss of generality, suppose the origin $0$ is the basepoint of
$\Gamma_i$ for all $i$.  By hypothesis, we know that $0
\in CH(\Lambda_{\Gamma_i})$ for all $i$, where
$CH$ denotes convex hull.  Then $T(\Gamma,R)$ contains all limits of rays from
$0$ to
$\Lambda_{\Gamma_i}$.  Therefore, $$\lim \sup \Lambda_{\Gamma_i} \subset
\overline{T(\Gamma,R)} \cap S^2_\infty \subset \Lambda_{\Gamma}.$$

Since $\Lambda_\Gamma \subset \lim \inf \Lambda_{\Gamma_i}$ for all $i$, the
result follows.
\end{proof}

From the proof of Theorem \ref{varylimset}, we can see that if for large
$i$, the basepoint of $\Gamma_i$ lies in the $K$-neighborhood of
$CH(\Lambda_{\Gamma_i})$, then $T(\Gamma, R+K)$ contains all limits of rays from
the basepoint of $\Gamma_i$ to $\Lambda_{\Gamma_i}$.  Then we know $$\lim \sup
\Lambda_{\Gamma_i}
\subset \overline{T(\Gamma,R+K)} \cap S^2_\infty \subset \Lambda_{\Gamma},$$
and hence $\Lambda_{\Gamma_i}$ converges to $\Lambda_\Gamma$ in the Hausdorff
topology.

Without loss of generality, suppose the origin $0$ is the basepoint of
$\Gamma_i$.  We will show that for large $i$, the origin $0$ lies within a
uniformly bounded neighborhood of $CH(\Lambda_{\Gamma_i})$.

Because $\Gamma$ is a limit of torsion-free Kleinian groups, $\Gamma$ itself
is torsion-free. (Lem 3.1.4, Canary-Epstein-Green \cite{can/ep/green})  Then
since $\Gamma$ is also nonabelian, we can conclude that $\Gamma$ is
nonelementary.  Therefore, $\Gamma$ contains a hyperbolic element $\gamma$. 
(Prop E.1, Maskit
\cite{maskit})  Let $A_\gamma$ denote the axis of $\gamma$ in
$\HHH$.  Then there exists a constant $K$ such that  
$d(0,A_\gamma) \leq K$.  As an element of the geometric limit, there exists a
sequence \{$\gamma_i \in \Gamma_i$\} such that $\gamma_i \rightarrow
\gamma$.  For large $i$, $\gamma_i$ is a hyperbolic element, and therefore, 
$A_{\gamma_i} \rightarrow A_\gamma$.  Then for large $i$, 
$d(0,A_{\gamma_i})
\leq K + 1$.  Because $A_{\gamma_i} \subset CH(\Lambda_{\Gamma_i})$ for all
$i$,  we can conclude that
$d(0,CH(\Lambda_{\Gamma_i})) \leq K + 1$ for large
$i$.  This concludes the proof of Corollary \ref{limitset2}.
\end{proof}

\bibliographystyle{amsplain}

\end{document}